\documentclass[12pt]{article}
\usepackage{amssymb}
\usepackage{amsfonts}
\usepackage{times}
\usepackage{cancel}
\usepackage{mathptmx}
\usepackage{amsmath}
\usepackage[usenames]{color}
\usepackage{mathrsfs}
\usepackage{amsfonts}
\usepackage{amssymb,amsmath}
\usepackage{amscd}
\usepackage{cite}
\usepackage{cases}
\usepackage{amsthm}
\def\cl{\centerline}
\def\la{\lambda}

\def\vs{\vspace*}
\def\Z{\mathbb{Z}}

\def\A{\mathcal{A}}
\def\L{\mathcal{L}}
\def\O{\mathcal{O}}

\def\C{\mathbb{C}}

\def\Ker{{\rm Ker}}

\def\vs{\vspace*}

\numberwithin{equation}{section}
\newtheorem{theo}{Theorem}[section]
\newtheorem{defi}[theo]{Definition}
\newtheorem{coro}[theo]{Corollary}
\newtheorem{lemm}[theo]{Lemma}
\newtheorem{prop}[theo]{Proposition}

\newtheorem{ex}[theo]{Example}

\newtheorem{remark}[theo]{Remark}

\allowdisplaybreaks[4]
\begin{document}
\begin{center}
{\bf\large $\mathcal{O}$-operators and Nijenhius operators of associative conformal algebras}
%\footnote{%Supported by National Natural Science
%Foundation grants of China (Grant No. 11301109).
%Corresponding author: lmyuan@hit.edu.cn}
\end{center}

\cl{Lamei Yuan}

\cl{\small School of Mathematics, Harbin Institute of Technology, Harbin
150001, China}
\cl{\small E-mail: lmyuan@hit.edu.cn}
\vs{6pt}

\vs{8pt}

\small\footnotesize
\parskip .005 truein
\baselineskip 3pt \lineskip 3pt
\noindent{\bf Abstract:} We study $\mathcal{O}$-operators of associative conformal algebras with respect to conformal bimodules. As natural generalizations of $\mathcal{O}$-operators and dendriform conformal algebras, we introduce the notions of twisted Rota-Baxter operators and conformal NS-algebras.
We show that twisted Rota-Baxter operators give rise to conformal NS-algebras, the same as $\mathcal{O}$-operators induce dendriform conformal algebras. And we introduce a conformal analog of associative Nijenhius operators and enumerate main properties. By using derived bracket construction of Kosmann-Schwarzbach and a method of Uchino, we obtain a graded Lie algebra whose Maurer-Cartan elements are given by $\mathcal{O}$-operators. This allows us to construct cohomology of $\mathcal{O}$-operators. This cohomology can be seen as the Hochschild cohomology of
 an associative conformal algebra with coefficients in a suitable conformal bimodule.

%In this paper, we study $\mathcal{O}$-operators of associative conformal algebras with respect to conformal bimodules. First, we introduce the notions of twisted Rota-Baxter operators and conformal NS-algebras, which are natural generalizations of  $\mathcal{O}$-operators and dendriform conformal algebras, respectively. We show that twisted Rota-Baxter operators give rise to conformal NS-algebras, the same as $\mathcal{O}$-operators induce dendriform conformal algebras.
%Second, we enumerate main properties of Nijenhius operators of associative conformal algebras. In particular, we present that a Nijenhuis operator gives rise to a whole hierarchy of Nijenhuis operators, associative conformal algebra structures and trivial deformations of associative conformal algebras. Third, by using derived bracket construction of Kosmann-Schwarzbach and a method of  Uchino, we construct a graded Lie algebra whose Maurer-Cartan elements are given by $\mathcal{O}$-operators. This allows us to construct cohomology for an $\mathcal{O}$-operator. This cohomology can be seen as the Hochschild cohomology of
% an associative conformal algebra with coefficients in a suitable conformal bimodule.

 \vs{10pt}
\noindent{\bf Key words:} Associative conformal algebras, $\mathcal{O}$-operators, Twisted Rota-Baxter operators, Nijenhius operators, Cohomology
 \vs{5pt}

\noindent{\bf Mathematics Subject Classification (2000):} 16D20, 16D70, 17A30, 17B55, 17B70.

\parskip .001 truein\baselineskip 6pt \lineskip 6pt

\section{Introduction}

A {\it Rota-Baxter operator} on an algebra $A$ over a field of characteristic zero is a linear map $T:A \rightarrow A$ satisfying
the following {\it Rota-Baxter equation}
\begin{align}
T(x)T(y)=T(T(x)y+xT(y)+qxy),
\end{align}
for all $x,\ y\in A$, where $q$ is an element of the ground field, called the {\it weight} of $T$. The pair $(A, T)$ is called a {\it Rota-Baxter
algebra of  weight $q$}. Rota-Baxter operators were first introduced by Baxter \cite{B} in the study of the fluctuation theory in probability. Later, they were further
developed by Rota \cite{Ro}, Atkinson \cite{A} and Cartier \cite{C} during the process of finding their interrelations with combinatorics. Also, these operators
were studied in integrable systems in the context of classical and modified Yang-Baxter equations
\cite{BD,STS}. In particular, it was established in \cite{B2,STS,Ku} that the Rota-Baxter equation on a Lie
algebra is precisely the operator form of the classical Yang-Baxter equation.

Rota-Baxter operators on associative algebras have been extensively studied. Interrelations between Rota-Baxter operators and associative analogues of the classical Yang-Baxter equation were studied in \cite{A3,BGN2}. As an associative analogue of Poisson structures on a manifold, Uchino \cite{U1} introduced the notion of a generalized Rota-Baxter operator, also known as $\O$-operator, which is a natural generalization of Rota-Baxter operators in the presence of bimodules. Let $(A,\ast)$ be an associative algebra and let $M$ be an $A$-bimodule. A linear map $T:M\rightarrow A$ is called a {\it generalized Rota-Baxter operator} (or {\it $\O$-operator}) on $A$ with respect to the bimodule $M$ if it satisfies
\begin{align}
T(m)\ast T(n)=T(m\cdot T(n)+T(m)\cdot n)
\end{align}
for all $m,n\in M$, where $\cdot$ means the bimodule action. Especially, when $M=A$ and $\cdot=\ast$, $T$
is reduced to a Rota-Baxter operator of weight $0$. Such an operator gave rise to a Loday's dendriform algebra structure \cite{Loday} on $M$ generalizing the fact from Rota-Baxter operators \cite{A2}. Therefore, $M$ inherits an associative structure as well. Further, Uchino \cite{U1} introduced the notion of a twisted Rota-Baxter operator in the context of associative algebras as an operator analog of twisted Poisson structures \cite{SW}. It turned out that twisted Rota-Baxter operators give rise to Nijenhuis (NS-)algebras, which were introduced in \cite{L}. By using the derived bracket construction, Uchino \cite{U2} constructed a differential graded Lie algebra associated to bi-graded Hochschild complex and proved that Rota-Baxter type operators are solutions of Maurer-Cartan equations. Following Uchino's method, Das \cite{Das} further constructed a cohomology of $\O$-operators, which can be seen as the Hochschild cohomology of an associative algebra with coefficients in a suitable bimodule. In the paper \cite{LBS}, the authors introduced and studied compatible $\O$-operators.
More precise study about Rota-Baxter operators, we refer readers to the book \cite{G}.

The study of Nijenhuis operators on Lie algebras dated back to late 1970s. It was discovered by Gel'fand and Dorfman  \cite{GD1,GD2} that Nijenhuis operators are closely related to Hamiltonian pairs. Interrelations between  Nijenhuis operators and deformations of Lie algebras were presented by Dorfman \cite{Dor}. Carinena and his coauthors \cite{CGM} introduced an associative version of classical Nijenhuis identity. Let $N:A\rightarrow A$ be a linear map on an associative algebra $(A,\ast)$. The operator $N$ is called an {\it associative Nijenhuis operator}, if it satisfies
\begin{eqnarray}
N(a)\ast N(b)=N\big(N(a)\ast b+a\ast N(b)-N(a\ast b)\big),
\end{eqnarray}
where $a,b\in A$. The deformed multiplication $a\times_N b:=N(a)\ast b+a\ast N(b)-N(a\ast b)$
is a new associative multiplication, which is compatible with the original one. In
this sense, an associative Nijenhuis operator induces a quantum bi-hamiltonian system (see \cite{CGM}).

An algebraic formalization of the properties of the operator product expansion (OPE) in two-dimensional conformal field theory \cite{BPZ} gave rise to a new class of algebraic systems, vertex operator algebras \cite{BO,F2}. The notion of a Lie conformal algebra encodes the singular part of the
OPE which is responsible for the commutator of two
chiral fields \cite{KAC}. Roughly speaking, Lie conformal algebras correspond to vertex algebras
by the same way as Lie algebras correspond to their associative enveloping algebras.

%Conformal algebras were introduced in \cite{KAC} as a useful tool for studying vertex algebras
%appeared in two-dimensional conformal field theory in mathematical physics
%\cite{BPZ}. The notion of a Lie conformal algebra encodes the singular part of the
%operator product expansion which is responsible for the commutator of two
%chiral fields. The relation of Lie conformal algebras
%and vertex algebras is similar to that of Lie algebras and their universal
%enveloping algebras.

The structure theory of finite (i.e., finitely
generated as $\C[\partial]$-modules) associative and Lie conformal algebras was developed in \cite{DK} and
later generalized in \cite{BDK} for pseudoalgebras over a wide class of cocommutative
Hopf algebras. From the algebraic point of view, the notions of conformal algebras \cite{DK}, their representations \cite{CK}
and cohomologies \cite{BKV,D,KK} are higher-level analogues of the ordinary notions in
the pseudo-tensor category \cite{BD2} associated with the polynomial Hopf algebra
(see \cite{BDK} for a detailed explanation).

Some features of the structure theory of conformal algebras (and their representations)
of infinite type were also considered in a series of works \cite{BKL1,BKL2,DK2,KAC2,Re1,Re2,Z1,Z2}. In this field,
one of the most urgent problems is to describe the structure of conformal
algebras with faithful irreducible representation of finite type (these algebras could be
of infinite type themselves). In \cite{BKL1,KAC2}, the conjectures on the structure of such algebras
(associative and Lie) were stated. The papers \cite{BKL1,DK,Z2} contain confirmations of these
conjectures under some additional conditions. Another problem is to classify simple and semisimple conformal algebras of linear
growth (i.e., of Gel'fand-Kirillov dimension one).
This problem was solved for finitely generated associative conformal algebras which
contain a unit \cite{Re1,Re2}, or at least an idempotent \cite{Z1,Z2}. The structure theory of associative conformal algebras with finite faithful representation similar to those examples of conformal algebras stated in these papers was developed in \cite{KO}. %The objects appearing from
%the consideration of conformal algebras with faithful irreducible representation of finite
%type are similar to those examples of conformal algebras stated in these papers.

In the recent paper \cite{HB1}, Hong and Bai developed a bialgebra theory for associative conformal algebras, which can be viewed as a conformal analogue of associative bialgebras \cite{A3} and also as an associative analogue of conformal bialgebras \cite{Li}. In particular, they introduced
the notions of $\O$-operators of associative conformal algebras and dendriform conformal
algebras to construct (antisymmetric) solutions of associative conformal Yang-Baxter
equation. In the present paper, we aim to extend the study of associative Rota-Baxter operators \cite{U1,U2, Das} and associative Nijenhius operators \cite{CGM} to the conformal case, and present more precise
properties of $\O$-operators and Nijenhius operators on associative conformal algebras. We hope that the present paper reveals further interesting interconnections and provides additional motivation to study
these operators.

The paper is arranged as follows. In Section 2, we recall the definitions of Lie and associative conformal algebras and their (bi-)modules. Also, we write down the constructions of Gerstenhaber's Lie bracket on the graded space of all multilinear maps over arbitrary vector spaces, and the derived bracket of Kosmann-Schwarzbach on a differential graded Lie algebra.  In addition, we gather some facts which will be used in this article.

In Section 3, we investigate some properties of $\mathcal{O}$-operators. Let $T:M\rightarrow \A$ be an $\mathcal{O}$-operator $T$ on an associative conformal algebra $\A$ with respect to a conformal $\A$-bimodule $M$. We show that the graph of $T$ is a subalgebra of the semi-direct product algebra $\A\oplus_0 M$ and its lift $\hat{T}$ is a Rota-Baxter operator of $\A\oplus_0 M$. It was known that $T$ could make $M$ into an associative conformal algebra $M_{ass}$. We will further prove that $T$ induces a conformal $M_{ass}$-bimodule structure on $\A$. And we prove that $T$ is also an $\mathcal{O}$-operator on the commutator Lie conformal algebra $\A^L$ with respect to the representation $(M,\rho)$. Finally, we show that compatible $\mathcal{O}$-operators give rise to compatible dendriform conformal algebra structures.

In Section 4, we construct a twisted version of Section 3. More specifically, we introduce the notion of a twisted Rota-Baxter operator, which is a generalization of $\O$-operators and characterized by a $2$-cocycle, and we construct a new algebraic structure, called conformal NS-algebra. We show that it is related to twisted Rota-Baxter operators in the same way that dendriform conformal algebras are related to $\O$-operators. Most of the results in Section 3 hold in the twisted case.

In Section 5, we introduce the notion of a Nijenhuis operator on associative conformal algebras and enumerate main properties.  First, we prove that Nijenhuis operators induce conformal NS-algebras. Second, we show that a Nijenhuis operator gives rise to a whole hierarchy of Nijenhuis operators and associative conformal algebra structures.
Third, we present interrelations between Nijenhuis operators and compatible
$\O$-operators, and interrelations between Nijenhuis operators of associative and Lie conformal algebras.  Finally, we show connections between Nijenhuis operators and deformations of associative conformal algebras.

In Section 6, we first recall Wu's construction of a differential graded Lie algebra structure on the Hochschild complex of an associative conformal algebra $\A$ by Gerstenhaber-bracket \cite{Ger}. Then we consider the semi-direct product algebra $(\A\oplus_0 M,\hat{\theta}_\la)$ of
$\A$ and a conformal $\A$-bimodule $M$, where $\hat{\theta}_\la$ is the associative $\la$-multiplication of $\A\oplus_0 M$.
The Hochschild complex $C^\bullet(\A\oplus_0 M)$ becomes a differential graded Lie algebra by Gerstenhaber-bracket and the coboundary map $d_{\hat{\theta}}:=[\hat{\theta},\cdot]$. Further, we define, due to \cite{KS1}, a derived bracket on $C^\bullet(\A\oplus_0 M)$ by
\begin{align*}
[[f,g]]:=(-1)^{{\rm deg}\,f}[[\hat\theta,f],g].
\end{align*}
Here the new bracket is a graded Lie bracket on $C^\bullet(M,\A)\subset C^\bullet(\A\oplus_0 M)$. We show that an element $T\in C^{1}(M, \A)$ is an $\O$-operator if and only if $T$ satisfies the Maurer-Cartan equation, i.e., $[[T,T]]=0.$ Also, an $\O$-operator
$T$ induces a differential $d_T:=[[T,\cdot]]$, which makes the graded Lie algebra
$(C^\bullet(M, \A), [[\cdot,\cdot]])$ into a differential graded Lie algebra. Hence we obtain a cohomology of the $\O$-operator $T:M\rightarrow \A$. This cohomology coincides with the Hochschild cohomology of $M$ with coefficients in $\A$.

Throughout this paper, all the vector spaces, linear maps and tensor products are over the complex field $\C$. Denote by $\mathbb{Z}$ the ring of integers and $\mathbb{N}$ the set of natural numbers. The elements of the vector space $\A$ are usually denoted by $a,b,c,\cdots$ and the elements of $M$ by $m,n,l,u,v,u_1,u_2,\cdots.$

\section{Preliminaries}

In this section, we recall some basic notions of associative and Lie conformal algebras along with their conformal modules and cohomology. We review  Gerstenhaber's construction of graded Lie algebra structure on the graded vector space of all multilinear maps over arbitrary vector spaces and derived bracket construction of Kosmann-Schwarzbach. Also, we gather some known results for later use. The material can be found in \cite{DK,D,Ger,KAC,KK,KS1,KS2,HB1}.
\subsection{Conformal algebras and modules}

\begin{defi}\label{defi1}
{\rm A { conformal algebra} $\mathcal{A}$ is a $\mathbb{C}[\partial]$-module  endowed with a $\mathbb{C}$-bilinear map
\begin{equation*}
\mathcal{A}\otimes \mathcal{A} \rightarrow \mathcal{A} [\lambda], ~~~~~~a\otimes b\mapsto a_\lambda b,
\end{equation*}
satisfying the following axiom
\begin{eqnarray}\label{sesquilinearity}
(\partial a)_\lambda b=-\lambda a_\lambda b,~~~
a_\lambda (\partial b)=(\partial+\lambda)a_\lambda b, ~~~~(\text{conformal~~sesquilinearity})
\end{eqnarray} for all $a, b\in \mathcal{A}$. If, in addition, it satisfies
\begin{eqnarray}\label{ASS}
(a_\lambda b)_{\lambda+\mu}c=a_\lambda (b_\mu c), ~~~~(\text{associativity})
\end{eqnarray}
for all $a, b\in \mathcal{A}$, then $\mathcal{A}$ is called an { associative conformal algebra}.
}
\end{defi}

\begin{defi}\rm
A { Lie conformal algebra} $\L$ is a $\C[\partial]$-module endowed with a $\C$-bilinear map
$$ \L\otimes \L\rightarrow \L[\lambda],\ \  a\otimes b \mapsto [a_\lambda b],$$
called the $\la$-bracket, and
satisfying the following axioms
\begin{align}
[\partial a_\lambda b]&=-\lambda[a_\lambda b],\ \ [ a_\lambda \partial b]=(\partial+\lambda)[a_\lambda b], \ \ \mbox{(conformal\  sesquilinearity)}\label{Lc1}\\
{[a_\lambda b]} &= -[b_{-\lambda-\partial}a], \ \ \mbox{(skew-symmetry)}\label{Lc2}\\
{[a_\lambda[b_\mu c]]}&=[[a_\lambda b]_{\lambda+\mu
}c]+[b_\mu[a_\lambda c]],\ \ \mbox{(Jacobi \ identity)}\label{Lc3}
\end{align}
for all $a, b, c\in \L$.
\end{defi}

Let $\mathcal{A}$ be an associative conformal algebra. It is well-known (see \cite{DK}) that the following $\la$-bracket
\begin{eqnarray}\label{L11}
[a_\la b]^L:=a_\lambda b-b_{-\la-\partial}a, \ \ \forall \ a,b\in\A
\end{eqnarray}
makes $\mathcal{A}$ into a Lie conformal algebra, which is called the {\it commutator} (or {\it sub-adjacent) Lie conformal algebra} of $\A$. We denote this Lie conformal algebra by $\A^L$.

Let $U$ and $V$ be two $\mathbb{C}[\partial]$-modules. We define the tensor product $U\otimes V$ of $\mathbb{C}[\partial]$-modules as the ordinary tensor product with
$\mathbb{C}[\partial]$-module structure ($u\in U,v \in V $):
\begin{eqnarray}
\partial(u\otimes v)=\partial u\otimes v+u\otimes \partial v.
\end{eqnarray}

\begin{defi}{\rm
Let $U$, $V$ and $W$ be $\mathbb{C}[\partial]$-modules.
\begin{itemize}
\item[(1)]A { left conformal linear map} from $U$ to $V$ is a $\mathbb{C}$-linear map
$f_\lambda: U\rightarrow V[\lambda]$, such that
\begin{eqnarray}
f_\lambda(\partial u)=-\lambda f_\lambda u, ~~ \forall ~~ u\in U.
\end{eqnarray}
\item[(2)]A { right conformal linear map} from $U$ to $V$ is a $\mathbb{C}$-linear map
$f_\lambda: U\rightarrow V[\lambda]$, such that
\begin{eqnarray}
f_\lambda(\partial u)=(\partial+\lambda)f_\lambda u, ~~ \forall ~~ u\in U.
\end{eqnarray}
A right conformal linear map is usually called a {conformal linear map} in short.
 \item[(3)]A { conformal bilinear map} from $U\otimes V$ to $W$ is a $\mathbb{C}$-bilinear map
$f_\lambda: U\otimes V\rightarrow W[\lambda]$, such that
\begin{eqnarray}
f_\lambda(\partial u,v)=-\lambda f_\lambda(u,v), ~~~
f_\lambda( u,\partial v)= (\partial+\lambda)f_\lambda(u,v),
\end{eqnarray}  for all $u\in U$ and $v\in V$.
\end{itemize}
}\end{defi}

Let $U$ and $V$ be two $\mathbb{C}[\partial]$-modules. We denote by $Chom(U,V)$ the vector space of all conformal linear maps from $U$ to $V$. It has a canonical structure of a $\C[\partial]$-module by
\begin{align}
(\partial f)_\la=-\la f_\la, \ \forall\ \ f_\la\in Chom(U,V).
\end{align}
In the special case $U=V$, we will write ${\rm Cend}(V)$ for $Chom(V, V)$. If $V$ is a finite $\C[\partial]$-module, then the $\C[\partial]$-module ${\rm Cend}(V)$ has a canonical structure of an associative conformal algebra defined by
\begin{align}
(f_\la g)_\mu v=f_\la(g_{\mu-\la}v),
\end{align}
for all $v\in V$ and $f,g\in {\rm Cend}(V)$. Further, the $\la$-bracket given by
\begin{align*}
[f_\la g]_\mu v=f_\la(g_{\mu-\la}v) -g_{\mu-\la}(f_\la v)
\end{align*}
defines a Lie conformal algebra structure on ${\rm Cend}(V)$. This is called the {\it general conformal algebra} on $V$ and denoted by $gc (V)$.

\begin{defi}{\rm  Let $\mathcal{A}$ be an associative conformal algebra and $M$ a $\C[\partial]$-module.
\begin{itemize}
\item[(1)] $M$ is called a left conformal module of $\A$ if the $\C$-bilinear map $\mathcal{A}\times M \rightarrow M[\lambda]$, $(a,m)\mapsto a_\lambda m$, is conformal sesquilinear and satisfies
\begin{eqnarray}\label{L}
(a_\lambda b)_{\lambda+\mu} m=a_\lambda(b_\mu m),
\end{eqnarray}
for all $a,b\in \mathcal{A}$ and $m\in M$.

\item[(2)] $M$ is called a right conformal module of $\A$ if the $\C$-bilinear map $\mathcal{A}\times M \rightarrow M[\lambda]$, $(m,a)\mapsto m_\lambda a$, is conformal sesquilinear and satisfies
\begin{eqnarray}\label{r}
(m_\lambda a)_{\lambda+\mu} b=m_\lambda(a_\mu b),
\end{eqnarray}
for all $a,b\in \mathcal{A}$ and $m\in M$.
\item[(3)] $M$ is called a conformal bimodule of $\A$ (or conformal $\A$-bimodule) if it is both a left conformal module and a right conformal module, and satisfies the following compatible condition
\begin{eqnarray}\label{B}
(a_\lambda m)_{\lambda+\mu} b=a_\lambda(m_\mu b),
\end{eqnarray}
for all $a,b\in \mathcal{A}$ and $m\in M$.
\end{itemize}}
\end{defi}

It follows that an associative conformal algebra $\mathcal{A}$ is a conformal bimodule over itself with the left and right $\la$-actions given by the $\la$-multiplication of $\mathcal{A}$. We call this conformal bimodule as a {\it adjoint bimodule}.

\begin{defi}\label{def1} \rm A  conformal module $V$ over a Lie conformal algebra $\L$
is a $\mathbb{C}[\partial]$-module endowed with a $\C$-bilinear map
$\L\otimes V\rightarrow V[\lambda]$, $(a,v)\mapsto a_\lambda v$, subject to the following conditions
\begin{eqnarray*}
&&(\partial a)_\lambda v =-\lambda a_\lambda v, \ \  a_\la (\partial v)=(\partial+\la)a_\lambda v,\\
&&[a_\lambda b]_{\lambda+\mu}v=a_\lambda (b_\mu v)-b_\mu (a_\lambda v),
\end{eqnarray*}
for all $a,b\in \L$ and $v\in V.$
\end{defi}
A conformal module $V$ over a Lie (or associative) conformal algebra $\L$ is called {\it finite}
if $V$ is finitely generated over $\mathbb{C}[\partial]$.
It is easy to see that a conformal module $V$ over a Lie conformal algebra $\L$ is the same as a homomorphism of Lie conformal algebras $\rho:\L\rightarrow gc(V)$, which is called a {\it representation} of $\L$ in the $\mathbb{C}[\partial]$-module $V$.

Let $\A$ be an associative conformal algebra. In the following, we consider a decomposition of $\A$ into a direct sum of two $\C[\partial]$-modules $\A_1$ and $\A_2$, namely, $\A=\A_1\oplus\A_2$ such that $\partial^{\A}=\partial^{\A_1}\oplus\partial^{\A_2}$.
The triple $(\A, \A_1,\A_2)$ is called a {\it matching pair} of
associative conformal algebras if $\A_1$ and $\A_2$ are subalgebras of $\A$ (cf.\cite{HB1,H1}). If an associative algebra decomposes into two subalgebras, it is also called an {\it associative twilled algebra} or simply {\it twilled algebra} in the literature (see, for example, \cite{CGM}). If a Lie algebra decomposes into two subalgebras, it is called a {\it twilled Lie algebra} in \cite{KS1}.

In the sequel, we denote the matching pair $(\A, \A_1,\A_2)$ of associative conformal algebras by $\A_1\bowtie\A_2$. One can easily check that $\A_1\bowtie\A_2$ is a matching pair of associative conformal algebras if and only if $\A_1$ (resp. $\A_2$)
is a conformal $\A_2$-bimodule (resp. $\A_1$-bimodule). In general, the associative $\la$-multiplication on $\A_1\bowtie\A_2$ has the form
\begin{eqnarray}\label{match}
(a,x)_\la (b,y)=(a_\la b+x\cdot^2_\la b+a\cdot^2_\la y, a\cdot^1_\la y+x\cdot^1_\la b+x_\la y),
\end{eqnarray}
where $a,b\in\A_1$, $x,y\in\A_2$, and the $\la$-action $\cdot^1_\la$ (resp. $\cdot^2_\la$) is the bimodule action of $\A_1$ on $\A_2$ (resp. $\A_2$ on $\A_1$).

The following is a special case of matching pairs of associative conformal algebras.
\begin{prop}\label{p1} (\cite{HB1}) Given an associative conformal algebra $\mathcal{A}$ and a conformal $\mathcal{A}$-bimodule $M$, the $\C[\partial]$-module $\mathcal{A}\oplus M$ carries an associative conformal algebra structure given by
\begin{align}\label{sum}
(a,m)_\lambda(b,n)=(a_\lambda b, a_\lambda n+m_\lambda b),
\end{align}
for all $a,b\in\mathcal{A}$ and $m,n\in M.$
\end{prop}

The associative conformal algebra from the proposition is called the {\it semi-direct product algebra} of $\A$ and $M$, and denoted by $\A\oplus _0 M$. It will be frequently used in this article.

\subsection{Cohomology of associative conformal algebras}
Let us describe the Hochschild cohomology complex $C^\bullet(\A,M)$ for an associative
conformal algebra $\A$ with coefficients in a conformal $\A$-bimodule $M$ by means of $\lambda $-products (see \cite{D,KK} for details). For any positive integer $n$, the space of $n$-cochains
$C^n(\A, M)$ consists of all multilinear maps of the form
\begin{align*}
\varphi_{\lambda_1,\cdots,\lambda_{n-1}}:\,\A^{\otimes n}&\longrightarrow M[\lambda_1, \cdots , \lambda_{n-1}]\\
a_1\otimes\cdots\otimes a_n&\longmapsto \varphi_{\lambda_1,\cdots,\lambda_{n-1}}(a_1,\cdots,a_n)
\end{align*}
satisfying the following sesquilinearity conditions:
\begin{align}
\varphi_{\lambda_1,\cdots,\lambda_{n-1}}(a_1,\cdots,\partial a_i,\cdots,a_n)&=-\la_i\varphi_{\lambda_1,\cdots,\lambda_{n-1}}(a_1,\cdots,a_n), ~~i=1,\cdots,n-1,\label{co1}\\
\varphi_{\lambda_1,\cdots,\lambda_{n-1}}(a_1,\cdots,\partial a_n)&=(\partial+\la_1+\cdots+\la_{n-1})\varphi_{\lambda_1,\cdots,\lambda_{n-1}}(a_1,\cdots,a_n).\label{co2}
\end{align}
The conformal Hochschild differential $\mathbf{d}:C^n(\A, M)\rightarrow C^{n+1}(\A, M)$ is defined by
\begin{align}\label{co3}
(\mathbf{d}\varphi)_{\lambda_1,\cdots,\lambda_{n}}(a_1,\cdots,a_{n+1})\,=\,&a_{1\,\la_1}\varphi_{\lambda_2,\cdots,\lambda_n}(a_2,\cdots,a_{n+1})\nonumber\\
&+\sum_{i=1}^{n}(-1)^i \varphi_{\lambda_1,\cdots,\lambda_i+\lambda_{i+1},\cdots,\lambda_{n}}(a_1,\cdots,a_{i\,{\la_i}\,}a_{i+1},\cdots,a_{n+1})\nonumber\\
&+(-1)^{n+1}\varphi_{\lambda_1,\cdots,\lambda_{n-1}}(a_1,\cdots,a_n)_{\la_1+\cdots+\la_n}a_{n+1}.
\end{align}
An $n$-cochain $\varphi\in C^n(\A, M)$ is called an {\it $n$-cocycle} if $\mathbf{d}\varphi=0$ and an element of the form $\mathbf{d}\varphi$, where $\varphi\in C^{n-1}(\A, M)$, is called an {\it $n$-coboundary}.
Denote by $Z^n(\A, M)$ and $B^n(\A, M)$ the subspaces of $n$-cocycles and $n$-coboundaries,
respectively. Then the quotient space $$H^n(\A, M) = Z^n(\A, M)/B^n
(\A, M)$$ is called the {\it $n$th Hochschild cohomology group} of $\A$ with coefficients in $M$.

For example, the space of $1$-cocycles $Z^1(\A, M) = \Ker \,{\mathbf{d}} \subseteq  C^1(\A, M)$ consists of all
$\C[\partial]$-linear maps $\varphi: \A \rightarrow M$ such that
\begin{align}\label{1-cocycle}
0=(\mathbf{d}\varphi)_{\lambda}(a,b)=a_\lambda \varphi (b)-\varphi(a_\lambda b)+\varphi(a)_{\lambda}b,
\end{align}
and the space of $2$-cocycles $Z^2(\A, M) = \Ker \,{\mathbf{d}} \subseteq  C^2(\A, M)$ consists of all conformal
sesquilinear maps $\varphi_{\lambda}: \A\otimes\A \rightarrow M[\la]$ such that
\begin{align}\label{2-cocycle}
0=(\mathbf{d}\varphi)_{\lambda,\mu}(a,b,c)=a_\lambda \varphi_\mu (b,c)-\varphi_{\lambda+\mu}(a_\lambda b,c)+\varphi_{\lambda}(a, b_\mu c)-\varphi_{\lambda}(a,b)_{\lambda+\mu}c.
\end{align}
\begin{remark}{\rm It is easy to see that $C^\bullet(\A,M)$ coincides with the reduced complex described in \cite{BKV}, where $C^n(\A,M)$ consists of adjacent classes of sesquilinear maps $\gamma_{\la_1,\cdots,\la_n}:\A^{\otimes n}\rightarrow M[\la_1,\cdots,\la_n]$ modulo the multiples of $(\partial+\la_1+\cdots+\la_n)$. The correspondence is given by
\begin{align*}
\gamma_{\la_1,\cdots,\la_n}\leftrightarrow \varphi_{{\la_1,\cdots,\la_{n-1}}}=\gamma_{\la_1,\cdots,\la_{n-1},-\partial-\la_1-\cdots-\la_{n-1}}.
\end{align*}}
\end{remark}

Recall that a conformal null extension of an associative conformal algebra $\A$ by means of a conformal bimodule $M$ over $\A$ is an associative conformal algebra $E$ in a short exact sequence
\begin{eqnarray*}
0\longrightarrow M \longrightarrow  E \longrightarrow \A \longrightarrow  0,
\end{eqnarray*}
such that $E$ is isomorphic to $\A\oplus M$ as a $\C[\partial]$-module and $M_\la M = 0$ in $E.$ Two conformal null extensions $E_1$ and $E_2$ are equivalent if there exists an isomorphism $E_1\rightarrow E_2$ such that the diagram
\begin{equation*}
\begin{CD}
0@>>> M @>>> E_1 @>>> \A @>>> 0\\
@. @V{{\rm id}_M}VV @VVV @V{{\rm id}_{\A}}VV\\
0@>>> M @>>> E_2 @>>> \A @>>> 0\\
\end{CD}
\end{equation*}
is commutative.

\begin{theo} \label{th10}(\cite{BKV,D}) Equivalence classes of conformal null
extensions of $\A$ by means of $M$ are in one-to-one correspondence with
 the elements of $H^2(\A, M)$.
\end{theo}

The semi-direct product algebra $\A\oplus_0 M$ appears as the trivial extension of $\A$ by $M$.
In general, given a conformal null extension of $\A$ by $M$, an associative $\lambda$-multiplication on $\A\oplus M$ has
the following form (cf.\cite{D,KK}):
\begin{align}\label{twist-sum}
(a,m)\circ^\varphi_\lambda(b,n)=\big(a_\lambda b, a_\lambda n+m_\lambda b+\varphi_\lambda(a,b)\big),
\end{align}
where $\varphi_\lambda$ is a $2$-cocycle in $C^2(\A, M)$. We denote the associative conformal
algebra $\A\oplus M$ equipped with the twisted $\lambda$-multiplication in \eqref{twist-sum} by $\A\oplus_{\varphi}M$, which will be studied in Section 4.

\subsection{Gerstenhaber-bracket and derived bracket}

\begin{defi}{\rm A differential graded Lie algebra (or simply dg-Lie algebra) is a triple  $(A,[\cdot,\cdot],d)$ such that
\begin{itemize}
\item[(1)] $A=\oplus_{i\in \mathbb{N}}A_i$, where $(A_i)_{i\in \mathbb{N}}$ is a family of $\C$-vector spaces, $[\cdot,\cdot]:A\times A\rightarrow A$ is a bilinear map of degree 0 and $d:A_k\rightarrow A_{k+1}$, $k\in \mathbb{N}$, is a graded
homomorphism of degree $+1$ such that $d^2=0$. An element $a\in A_k$
is said to be homogeneous of degree $k={\rm deg}\,a$.
\item[(2)] $[\cdot,\cdot]: A\times A\rightarrow A$ defines a structure of graded Lie algebra, i.e., for homogeneous elements $a,b,c\in A$ there hold
    (I) Graded commutativity:
    \begin{eqnarray*}
     [a,b]=-(-1)^{{\rm deg\,}a\, {\rm deg\,}b}[b,a],
\end{eqnarray*}
(II) Graded Jacobi identity:
\begin{align*}
(-1)^{{\rm deg\,}a\, {\rm deg\,}c}[[a,b],c]+(-1)^{{\rm deg\,}b\, {\rm deg\,}a}[[b,c],a]+(-1)^{{\rm deg\,}c\, {\rm deg\,}b}[[c,a],b]=0.
\end{align*}
\item[(3)]  $d$ is compatible with the graded Lie algebra structure, i.e.,
\begin{align*}
d([a,b])=[d(a),b]+(-1)^{{\rm deg}\,a}[a,d(b)].
\end{align*}
 \end{itemize}
}\end{defi}
The above graded Jacobi identity is equivalent to
\begin{eqnarray}\label{Leibniz}
    [a,[b,c]]=[[a,b],c]+(-1)^{{\rm deg\,}a\, {\rm deg\,}b}[b,[a,c]],
\end{eqnarray}
which is called {\it graded Leibniz identity}, sometimes also called {\it graded Loday identity}.

\begin{defi}{\rm Let $(A,[\cdot,\cdot],d)$ be a dg-Lie algebra and $a\in A_1$. We say that $a$ is a Maurer-Cartan element in $(A,[\cdot,\cdot],d)$ if it verifies the Maurer-Cartan equation, i.e.,
\begin{align}
d(a)+\frac12[a,a]=0.
\end{align}
Further, we say that $a$ is a strong Maurer-Cartan element in $(A,[\cdot,\cdot],d)$ if it satisfies
\begin{align}
d(a)=\frac12[a,a]=0.
\end{align}
}\end{defi}

Let $(A,[\cdot,\cdot],d)$ be a dg-Lie algebra. Define a new bracket on $A$ by
\begin{align*}
[a,b]_d=(-1)^{{\rm deg}(a)}[d(a),b], \ \forall \ a,b\in A,
\end{align*}
called the {\it derived bracket}.
Then the new bracket becomes a Leibniz bracket (or  Loday bracket), namely, it satisfies the graded Leibniz identity \eqref{Leibniz}. This method of constructing a new product is called a {\it derived bracket construction} of Kosmann-Schwarzbach (\cite{KS1,KS2}). The derived bracket construction plays important
roles in modern analytical mechanics and Poisson geometry. It is known that
several important brackets, e.g., Poisson brackets, Schouten-Nijenhuis brackets, Lie
algebroid brackets, Courant brackets and BV-brackets are induced by the derived
bracket construction.

%The new bracket is called a {\bf derived bracket} (cf. \cite{KS1,KS2}, see also \cite{U3}). It is well-known that the derived bracket satisfies the graded Leibniz identity.

The following basic lemma given in \cite{KS1} will be used in Section 6.
\begin{lemm}\label{lemm55}Let $(A,[\cdot,\cdot],d)$ be a dg-Lie algebra and let $\mathfrak{h}\subset A$ be an abelian subalgebra, i.e., $[\mathfrak{h},\mathfrak{h}]=0$. If the derived bracket $[\cdot,\cdot]_d$ is closed in $\mathfrak{h}$, then $(\mathfrak{h}, [\cdot,\cdot]_d)$ forms a graded Lie algebra.
\end{lemm}

Finally, let us recall Gerstenhaber's construction of a graded Lie algebra structure on the graded vector space of all multilinear maps on an arbitrary vector space $V$ (cf.\cite{Ger}). For $n\in \mathbb{N}$, set $\mathfrak{g}^n(V):=Hom(V^{\otimes n}, V)$, which consists of all $n$-linear maps. If $f\in \mathfrak{g}^n(V)$, then $f$ is said to be homogenous element of degree $n$.  Define a graded Lie bracket which is called the {\it Gerstenhaber-bracket}, or shortly {\it G-bracket} on $\mathfrak{g}^\bullet(V):=\oplus_{n\in \mathbb{N}} \mathfrak{g}^n(V)$ by
\begin{align*}
[f,g]=f\circ g-(-1)^{(m-1)(n-1)} g\circ f,
\end{align*}
where $\circ$ is the composition of maps defined by
\begin{align*}
(f\circ g)(v_1,\cdots,v_{m+n-1})=\sum_{i=1}^{m}(-1)^{(i-1)(n-1)}f(v_1,\cdots,v_{i-1},g(v_i,\cdots,v_{i+n-1}),v_{i+n},\cdots,v_{m+n-1}),
\end{align*}
for all $f\in\mathfrak{g}^m(V)$, $g\in\mathfrak{g}^n(V)$ and $v_1,\cdots,v_{m+n-1}\in V.$ Notice that the degree of $f\circ g$ is $m+n-1$. That is to say the G-bracket is of degree $-1$.
There hold two fundamental identities:

(i) Graded commutativity:
\begin{align*}
[f,g]=-(-1)^{(m-1)(n-1)}[g,f];
\end{align*}

(ii) Graded Jacobi identity:
\begin{align*}
(-1)^{(m-1)(l-1)}[[f,g],h]+(-1)^{(l-1)(n-1)}[[h,f],g]+(-1)^{(n-1)(m-1)}[[g,h],f]=0,
\end{align*}
which is equivalent to the following graded Leibniz identity
\begin{align}\label{1-3}
[f,[g,h]=[[f,g],h]+(-1)^{(n-1)(m-1)}[g,[f,h]],
\end{align}
where $f\in\mathfrak{g}^m(V)$, $g\in\mathfrak{g}^n(V)$ and $h\in \mathfrak{g}^l(V)$.

\section{$\mathcal{O}$-operators and dendriform conformal algebras}
\vs{8pt}

In this section,
we investigate some properties of $\mathcal{O}$-operators on associative conformal algebras with respect to conformal bimodules    and their connections with derivations and dendriform (Lie and left-symmetric) conformal algebras. We also introduce the notions of compatible $\mathcal{O}$-operators and compatible dendriform conformal algebras and describe their interrelations.

Let's start with recalling the definition of $\mathcal{O}$-operators given in \cite{HB1}.
\begin{defi} {\rm  Let $M$ be a conformal bimodule over an associative conformal algebra $\A$. A $\C[\partial]$-module homomorphism $T:M\rightarrow \A$ is called an {$\mathcal{O}$-operator} on $\A$ with respect to $M$ if it satisfies
\begin{eqnarray}\label{O}
T(m)_\lambda T(n)=T\big(T(m)_\lambda n+m_\lambda T(n)\big),
\end{eqnarray} for all $m,n\in M.$
}\end{defi}

When taking $M=\mathcal{A}$, an $\mathcal{O}$-operator $T$ is nothing but a {\it Rota-Baxter operator} on $\mathcal{A}$ (cf. \cite{HB1}), namely, $T$ satisfies
\begin{eqnarray}
T(a)_\lambda T(b)=T\big(T(a)_\lambda b+a_\lambda T(b)\big),
\end{eqnarray}
for all $a,b\in \A.$ Hence an  $\mathcal{O}$-operator $T$ on an associative conformal algebra $\A$ with respect to a conformal bimodule $M$ is also called a {\it generalized
Rota-Baxter operator}.

Let $M$ be a conformal bimodule over an associative conformal algebra $\A$.  By Proposition \ref{p1}, $\A\oplus_0 M$ is an associative conformal algebra with respect to \eqref{sum}. Assume that $T:M\rightarrow \A$ is a $\C[\partial]$-module homomorphism. We denote the graph of $T$ by $Gr(T)$,
\begin{align}
Gr(T)=\{(T(m),m)|m\in M\}.
\end{align}
%Then we have the following result.
\begin{prop}\label{p4} $T:M\rightarrow \A$ is an  $\mathcal{O}$-operator if and only if $Gr(T)$ is a subalgebra of
 $\A\oplus_0 M$.
\end{prop}
\begin{proof} For any $(T(m),m)$, $(T(n),n)\in Gr(T)$, we have
\begin{align*}
(T(m),m)_\lambda (T(n),n)=(T(m)_\lambda T(n), T(m)_\lambda n+m_\lambda T(n)).
%&=(T(T(m)_\lambda n+m_\lambda T(n)),T(m)_\lambda n+m_\lambda T(n))\\
%&=(T(k),k)\in Gr(T)[\lambda].
\end{align*}
Hence $T$ is an $\mathcal{O}$-operator if and only if $(T(m)_\lambda T(n), T(m)_\lambda n+m_\lambda T(n))$ is in $Gr(T)[\lambda]$.
 \end{proof}

Note that $Gr(T)$ and $M$ are isomorphic as conformal $\A$-bimodules by identification $(T(m),m)\cong m$. Hence, if $T$ is an $\mathcal{O}$-operator, i.e., $Gr(T)$ is an associative conformal subalgebra of $\A\oplus_0 M$, then $M$ is also an associative conformal algebra.

Given an arbitrary $\C[\partial]$-module homomorphism $T:M\rightarrow \A$, we define a lift of $T$, $\hat{T}$, as an endomorphism on $\A\oplus M$ by $\hat{T}(a,m):=(T(m),0)$, for all $a\in\A$ and $m \in M.$

\begin{prop} \label{p3-3} $T:M\rightarrow \A$ is an  $\mathcal{O}$-operator if and only if $\hat{T}$ is a Rota-Baxter operator (of weight 0) on  $\A\oplus_0 M$.
\end{prop}
\begin{proof} $\hat{T}$ is obviously a $\C[\partial]$-module homomorphism. For any $(a,m),~(b,n)\in \A\oplus M$, we have
\begin{align}\label{1-1}
\hat{T}(a,m)_\lambda \hat{T}(b,n)=(T(m),0)_\lambda (T(n),0)=(T(m)_\lambda T(n),0),
\end{align}
and
\begin{align}\label{1-2}
\hat{T}\big(\hat{T}(a,m)_\lambda (b,n)+(a,m)_\lambda \hat{T}(b,n)\big)&=\hat{T}\big((T(m),0)_\lambda (b,n)+(a,m)_\lambda (T(n),0)\big)\nonumber\\
&=\hat{T}\big((T(m)_\lambda b,T(m)_\lambda n)+(a_\lambda T(n),m_\lambda T(n))\big)\nonumber\\
&=\big(T(T(m)_\lambda n+m_\lambda T(n)),0\big).
\end{align}
Combining \eqref{1-1} with \eqref{1-2}, we obtain the result.
\end{proof}

Now, we recall the notion of a dendriform conformal algebra introduced in \cite{HB1}. It is a conformal analog of the classical dendriform algebras, which were first introduced by Loday \cite{Loday} with motivation from algebraic $K$-theory.
\begin{defi}{\rm A { dendriform conformal algebra} is a triple $(E,\succ_\lambda,\prec_\lambda)$ consisting of a $\C[\partial]$-module $E$ and two $\la$-multiplications $\succ_\lambda, \prec_\lambda:E\times E\rightarrow E[\lambda]$, which are conformal sesquilinear maps and satisfy the following axioms:
\begin{align}
a\succ_\lambda(b\succ_\mu c)&=(a\succ_\lambda b+a\prec_\lambda b)\succ_{\lambda+\mu}c,\label{cd1}\\
(a\prec_\lambda b)\prec_{\lambda+\mu}c&=a\prec_\lambda (b\succ_\mu c+b\prec_\mu c),\label{cd2}\\
(a\succ_\lambda b)\prec_{\lambda+\mu}c&=a\succ_\lambda (b\prec_\mu c),\label{cd3}
\end{align} for all $a,b,c\in E.$}
\end{defi}

It was shown in \cite{Loday} that given a dendriform algebra $(E,\succ,\prec)$, the sum of the two multiplications $$x\star y:=x\succ y+x\prec y$$ is associative. In the conformal case, the same holds.

\begin{prop}\label{p2} (\cite{HB1}) If $(E,\succ_\lambda,\prec_\lambda)$ is a dendriform conformal algebra, then  $(E,\star_\lambda)$ forms an associative conformal algebra, where $\star_\lambda$ is defined by
\begin{eqnarray}
a\star_\lambda b:=a\succ_\lambda b+a\prec_\lambda b, %~~ \forall ~~a,b\in E.
\end{eqnarray} for all $a,b\in E.$ The algebra $(E,\star_\lambda)$ is called the {\it associated associative conformal algebra} of $(E,\succ_\lambda,\prec_\lambda)$, and denoted by $E_{ass}$.
\end{prop}

%The algebra $(E,\star_\lambda)$ from Proposition \ref{p2} is called the {\bf associated associative conformal algebra} of $(E,\succ_\lambda,\prec_\lambda)$, and denoted by $E_{ass}$.

The following proposition says that an $\mathcal{O}$-operator has an underlying dendriform structure.

\begin{prop}  (\cite{HB1}) \label{p3} Let $T:M\rightarrow \A$ be an $\mathcal{O}$-operator. Then $M$ becomes a dendriform conformal algebra with the $\lambda$-multiplications given by
\begin{align}
m\succ^T_\lambda n=T(m)_\lambda n, ~~ m\prec^T_\lambda n=m_\lambda T(n),
\end{align}
where $m,n\in M.$ And there is an induced dendriform conformal algebra structure on $T(M)= \{T(m) | m \in M\} \subset\A$ given
by
\begin{align}
T(m)\succ_\lambda T(n)=T(m\succ^T_\lambda n), ~~ T(m)\prec_\lambda T(n)=T(m\prec^T_\lambda n), \ \forall \ m,n\in M.
\end{align}
If, in addition, $T$ is invertible, then there exists a dendriform conformal algebra structure on $\A$
defined by
\begin{align}
a\succ_\lambda b= T(a_\la T^{-1}(b)),~~  a \prec_\lambda b = T(T^{-1}(a)_\la b),
\end{align}
for all $a,b\in\A.$
\end{prop}

It follows from Propositions \ref{p2} and \ref{p3} that if $T:M\rightarrow \A$ is an $\mathcal{O}$-operator, then $M$ has an associative $\lambda$-product of the form
\begin{align}\label{M-ass}
m\star_\lambda n=T(m)_\lambda n+m_\lambda T(n), ~~\forall~m,n\in M.
\end{align}
We denote the associative conformal algebra $(M,\star_\lambda)$ by $M_{ass}$. Notice that \eqref{M-ass} implies that
\begin{align}\label{M-ass-1}
T(m\star_\lambda n)=T(m)_\lambda T(n), ~~\forall~m,n\in M.
\end{align}
Hence $T$ is an algebra homomorphism from $M_{ass}$ to $\A$.

\begin{lemm} \label{lemm1} Under the assumptions above, $\A$ becomes a conformal $M_{ass}$-bimodule
by the following $\lambda$-actions:
\begin{align}\label{2-1}
m\cdot_\lambda a=T(m)_\lambda a-T(m_\lambda a), ~~ a\cdot_\lambda m=a_\lambda T(m)-T(a_\lambda m),
\end{align}
where $m\in M_{ass}$ and $a\in\A.$
\end{lemm}
\begin{proof} It is easy to see that the two $\lambda$-actions defined by \eqref{2-1} are conformal sesquilinear maps. For any $m,n\in M$ and $a,b\in\A,$ we have
\begin{align*}
m\cdot_\lambda (n\cdot_\mu a)&=m\cdot_\lambda\big(T(n)_\mu a-T(n_\mu a)\big)\\
&=T(m)_\lambda(T(n)_\mu a)-T(m_\lambda (T(n)_\mu a))-T(m)_\lambda T(n_\mu a)+T(m_\lambda T(n_\mu a))\\
&=T(m)_\lambda(T(n)_\mu a)-T(m_\lambda (T(n)_\mu a))-\cancel{T(m_\lambda T(n_\mu a))}-T(T(m)_\lambda (n_\mu a))+\cancel{T(m_\lambda T(n_\mu a))}\\
&=T(m)_\lambda(T(n)_\mu a)-T(m_\lambda (T(n)_\mu a))-T(T(m)_\lambda (n_\mu a)).
\end{align*}
On the other hand, we have
\begin{align*}
(m\star_\lambda n)\cdot_{\lambda+\mu} a&=T(m\star_\lambda n)_{\lambda+\mu} a-T((m\star_\lambda n)_{\lambda+\mu} a)\\
&=(T(m)_\lambda T(n))_{\lambda+\mu} a-T((T(m)_\lambda n+m_\lambda T(n))_{\lambda+\mu} a)\\
&=T(m)_\lambda(T(n)_\mu a)-T(m_\lambda (T(n)_\mu a))-T(T(m)_\lambda (n_\mu a)),
\end{align*}
where we have used \eqref{M-ass} and \eqref{M-ass-1}. This proves that  $m\cdot_\lambda (n\cdot_\mu a)=(m\star_\lambda n)\cdot_{\lambda+\mu} a$.

 Similarly, we can obtain $(a\cdot_\lambda m)\cdot_{\lambda+\mu} n=a\cdot_\lambda (m\star_{\mu}n)$ and $m\cdot_\lambda (a\cdot_\mu n)=(m\cdot_\lambda a)\cdot_{\lambda+\mu} n$. Hence $\A$ is a conformal $M_{ass}$-bimodule.
\end{proof}

It follows from Lemma \ref{lemm1} that if  $T:M\rightarrow \A$ is an $\mathcal{O}$-operator, then we have a matching pair $\A\bowtie M_{ass}$ of associative conformal algebras. The associative $\la$-multiplication of $\A\bowtie M_{ass}$ has the form
\begin{align}\label{M-P}
(a,m)_\lambda(b,n)=(a_\lambda b+a\cdot_\la n+m\cdot_\la b, a_\lambda n+m_\lambda b+m\star_\la n),
\end{align}
for all $a,b\in\mathcal{A}$ and $m,n\in M_{ass},$ where $\cdot_\la$ means the conformal bimodule action of $M_{ass}$ on $\A$ defined by \eqref{2-1} and $\star_\la$ is the associative $\la$-multiplication of $M_{ass}$ defined by \eqref{M-ass}.

Recall that a $\C[\partial]$-linear map $d:\A\rightarrow M$ is called a {\it derivation} from an associative conformal algebra $\A$ to its conformal bimodule $M$ if it satisfies
\begin{align}\label{der}
d(a_\la b)=d(a)_\la b+a_\la d(b), \ \forall \ a,b\in\A.
\end{align}
It follows from \eqref{1-cocycle} that a derivation is exactly a $1$-cocycle in $C^1(\A,M)$.

The following proposition describes a close relation between $\mathcal{O}$-operators and derivations.

\begin{prop}  Let $T:M\rightarrow \A$ be an $\mathcal{O}$-operator and $\Omega:\A\rightarrow M$ a derivation satisfying
\begin{align}\label{der4}
\Omega(a)\star_\la \Omega(b)=\Omega\big(\Omega(a)\cdot_\la b+a\cdot_\la \Omega(b)\big),\ \ \forall\ a,b\in\A.
\end{align}
Then\begin{itemize}
\item[(1)] the composition map $T\Omega:\A\rightarrow \A$ satisfies
\begin{align}\label{der2}
T\Omega(a)_\la T\Omega(b)=T\Omega\big(T\Omega(a)_\la b+a_\la T\Omega(b)-T\Omega(a_\la b) \big),
\end{align}
for all $a,b\in\A.$
\item[(2)] The composition map $T\Omega T:M\rightarrow \A$ is a second $\mathcal{O}$-operator.
\end{itemize}
\end{prop}
\begin{proof} (1) By Lemma \ref{lemm1}, we have
\begin{align}\label{der3}
\Omega(a)\cdot_\la b+a\cdot_\la \Omega(b)&=T\Omega(a)_\la b-T(\Omega(a)_\la b)+a_\la T\Omega(b)-T(a_\la \Omega(b))\nonumber\\
&=T\Omega(a)_\la b+a_\la T\Omega(b)-T\Omega(a_\la b),
\end{align}
for all $a,b\in\A.$ Here the derivation condition of $\Omega$ is used. Applying $T\Omega$ to the both sides of \eqref{der3}, the
left-hand side reads
\begin{align*}
T\Omega\big(\Omega(a)\cdot_\la b+a\cdot_\la \Omega(b)\big)\stackrel{\eqref{der4}}{=}T(\Omega(a)\star_\la \Omega(b))\stackrel{\eqref{M-ass-1}}{=}T\Omega(a)_\la T\Omega(b),
\end{align*}
whereas the right-hand side obviously reads $T\Omega\big(T\Omega(a)_\la b+a_\la T\Omega(b)-T\Omega(a_\la b)\big)$. Hence we obtain \eqref{der2}.

(2) Put $a:=T(m)$ and $b:=T(n)$ for any $m,n\in M.$ Plugging this into \eqref{der2} gives
\begin{align}\label{3-6}
T\Omega T(m)_\la T\Omega T(n)=T\Omega\big(T\Omega T(m)_\la T(n)+T(m)_\la T\Omega T(n)-T\Omega(T(m)_\la T(n))\big).
\end{align}
As $T:M\rightarrow \A$ is an $\mathcal{O}$-operator, we have
\begin{align*}
T\Omega T(m)_\la T(n)&=T\big(T\Omega T(m)_\la n+\Omega T(m)_\la T(n)\big),\\
T(m)_\la T\Omega T(n)&=T\big(T(m)_\la \Omega T(n)+m_\la T\Omega T(n)\big),
\end{align*}
and as $\Omega:\A\rightarrow M$ is a derivation, we have
\begin{align*}
T\Omega(T(m)_\la T(n))=T\big(\Omega T(m)_\la T(n)+T(m)_\la \Omega T(n)\big).
\end{align*}
Then \eqref{3-6} becomes
\begin{align*}
T\Omega T(m)_\la T\Omega T(n)=&\,T\Omega T\big(T\Omega T(m)_\la n+\cancel{\Omega T(m)_\la T(n)}+\cancel{T(m)_\la \Omega T(n)}+m_\la T\Omega T(n)\big)\\&-T\Omega T\big(\cancel{\Omega T(m)_\la T(n)}+\cancel{T(m)_\la \Omega T(n)}\big)\\
=&\,T\Omega T\big(T\Omega T(m)_\la n+m_\la T\Omega T(n)\big),
\end{align*}
as required.
\end{proof}

\begin{remark}{\rm Condition \eqref{der4} is equivalent to that $\Omega:\A\rightarrow M_{ass}$ is an $\O$-operator. Condition \eqref{der2} actually says that $T\Omega$ is a Nijenhuis operator on $\A$ (see \eqref{N}). This implies a close interrelation between $\O$-operators and Nijenhuis operators,  which will be studied in Section 5.}
\end{remark}

\begin{ex} {\rm Given a dendriform conformal algebra $(E,\succ_\lambda,\prec_\lambda)$, we have an associative conformal algebra $E_{ass}$ by Proposition \ref{p2}. One can easily check that $E$ is a conformal $E_{ass}$-bimodule by
\begin{align}
e\cdot_\lambda x := e\succ_\lambda x, ~~ x \cdot_\lambda e = x \prec_\lambda e,
\end{align}
where $e\in E_{ass}$ and $x\in E.$ Under this setting, the identity map
${\rm {id} }: E\rightarrow E_{ass}$ is an $\mathcal{O}$-operator and the corresponding dendriform conformal
algebra is the original one, i.e., $(E,\succ_\lambda,\prec_\lambda)$. Hence all dendriform conformal algebras are induced by $\mathcal{O}$-operators.
}\end{ex}

\begin{defi}{\rm (\cite{HB2}) Let $\L$ be a Lie conformal algebra and $\rho:\L\rightarrow gc(V)$ a representation. If a $\C[\partial]$-module homomorphism $T:V\rightarrow \L$ satisfies
\begin{eqnarray}
[T(u)_\lambda T(v)]=T\big(\rho(T(u))_\lambda v-\rho(T(v))_{-\la-\partial}u\big),
\end{eqnarray}
for all $u,v\in V$, then $T$ is called an $\mathcal{O}$-operator of $\L$ associated with $\rho$.}
\end{defi}

 Let $\A$ be an associative
conformal algebra and $M$ a conformal bimodule over $\A$. We consider the
commutator Lie conformal algebra of $\A$, $\A^L$, which is defined by
\eqref{L11}. Then $M$ can be given a structure of Lie conformal algebra representation $\rho:\A^L\rightarrow gc(M)$ by
\begin{eqnarray}\label{L2}
\rho(a)_\lambda m:=a_\lambda m-m_{-\la-\partial}a,
\end{eqnarray}
where $a\in\A^L$ and $m\in M$. We denote this representation by $(M,\rho)$.

\begin{theo} \label{th2} Let $T:M\rightarrow \A$ be an $\mathcal{O}$-operator on an associative conformal algebra $\A$ with respect to a conformal bimodule $M$. Then $T$ is also an $\mathcal{O}$-operator on the commutator Lie conformal algebra $\A^L$ with respect to the representation $(M,\rho)$.
\end{theo}
\begin{proof}For any $m,n\in M$, we have
\begin{align*}
[T(m)_\lambda T(n)]&\stackrel{\eqref{L11}}{=}T(m)_\la T(n)-T(n)_{-\la-\partial}T(m)\\
&\stackrel{\eqref{O}}{=}T(T(m)_\la n+m_\la T(n))-T(T(n)_{-\la-\partial}m+n_{-\la-\partial}T(m))\\
&\,\,=\,T(T(m)_\la n-n_{-\la-\partial}T(m))-T(T(n)_{-\la-\partial}m-m_\la T(n))\\
&\stackrel{\eqref{L2}}{=}T(\rho(T(m))_\la n)-T(\rho(T(n))_{-\la-\partial}m).
\end{align*}
The proof is finished.
\end{proof}

Let $T$ be an $\mathcal{O}$-operator on an associative conformal algebra $\A$ with respect to a conformal $\A$-bimodule $M$. Then $M$ carries an associative conformal algebra structure $M_{ass}$ given by \eqref{M-ass} and there is a conformal $M_{ass}$-bimodule structure on $\A$ by Lemma \ref{lemm1}. 
Note that the commutator Lie conformal algebra structure $M_{ass}^{L}$ on $M_{ass}$ is given by
\begin{align}\label{untwist}
[m_\lambda n]^L\stackrel{\eqref{L11}}{=}m\star_\lambda n-n\star_{-\lambda-\partial}m\stackrel{\eqref{M-ass}}{=}T(m)_\lambda n+m_\lambda T(n)-T(n)_{-\lambda-\partial} m-n_{-\lambda-\partial} T(m),
\end{align}
and its Lie conformal algebra representation on $\A$ is given by $\rho_\A:M_{ass}^L\rightarrow gc(\A)$, where
\begin{align*}
\rho_{\A}(m)_\la a&\stackrel{\eqref{L2}}{=}m\cdot_\la a-a\cdot_{-\la-\partial}m\\
&\stackrel{\eqref{2-1}}{=}T(m)_\lambda a-T(m_\lambda a)-a_{-\la-\partial} T(m)+T(a_{-\la-\partial} m),
\end{align*}
for all $a\in\A$ and $m\in M.$

On the other hand, it follows from  Theorem \ref{th2} that $T$ induces an $\mathcal{O}$-operator on the commutator Lie conformal algebra $\A^L$ with respect to the representation $(M, \rho)$. Moreover, it is easy to check that the following $\la$-bracket
\begin{align*}
[m_\la n]:&=\rho(T(m))_\la n-\rho(T(n))_{-\la-\partial}m\\
&\stackrel{\eqref{L2}}{=}(T(m)_\la n-n_{-\la-\partial}T(m))-(T(n)_{-\la-\partial}m-m_{\la}T(n))
\end{align*}
makes $M$ into a Lie conformal algebra, denoted by $M^L_\rho$. Define $\rho_{\A}^{\prime}:M^L_\rho\rightarrow gc (\A)$ by
\begin{align*}
\rho_{\A}^{\prime}(m)_\la(a):&=[T(m)_\la a]^L+T(\rho(a)_{-\la-\partial} m)\\
\text(by~ \eqref{L11} ~ and ~ \eqref{L2}) &=T(m)_\la a- a_{-\la-\partial}T(m)+T(a_{-\la-\partial} m-m_\la a),
\end{align*}
where $a\in\A$ and $m\in M.$
It is not difficult to show that $\rho_{\A}^{\prime}$ is a representation of $M^L_\rho$ in $\A$.
We see that the two Lie conformal algebra structures $M^L_{ass}$ and $M^L_\rho$ are exactly the same, and the corresponding representations $\rho_{\A}$ and $\rho_{\A}^{\prime}$ on $\A$ are also the same.

Let's recall the definition of left-symmetric conformal algebras introduced in \cite{HB2}.

\begin{defi} {\rm A left-symmetric conformal algebra $\A$ is a $\C[\partial]$-module with a $\C$-bilinear map $\circ_\la:\A\times\A\rightarrow\A[\la], (a,b)\mapsto a\circ_\la b$, which is a conformal sesquilinear map and satisfies
\begin{align}
(a\circ_\la b)\circ_{\la+\mu} c-a\circ_\la(b\circ_\mu c)=(b\circ_\mu a)\circ_{\la+\mu}c-b\circ_\mu(a\circ_\la c),
\end{align}
for all $a,b,c\in\A.$}
\end{defi}

The following proposition says that dendriform conformal algebras give rise to left-symmetric structures.

\begin{prop}\label{P10} If $(E,\succ_\lambda,\prec_\lambda)$ is a dendriform conformal algebra, then $(E,\circ_\lambda)$ forms a left-symmetric conformal algebra, where $\circ_\lambda$ is defined by
\begin{eqnarray}
a\circ_\lambda b:=a\succ_\lambda b-b\prec_{-\la-\partial} a, ~~ \forall ~~a,b\in E.
\end{eqnarray}
\end{prop}
\begin{proof} Straightforward.
\end{proof}

 It follows from Propositions \ref{p3} and \ref{P10} that if $T:M\rightarrow\A$ is an $\O$-operator on an associative conformal algebra $\A$ with respect to a conformal bimodule $M$,  then $M$ carries a left-symmetric conformal algebra structure with the $\la$-multiplication defined by
\begin{eqnarray}
m\circ_\lambda n=T(m)_\lambda n-n_{-\la-\partial} T(m), ~~ \forall ~~m,n\in M.
\end{eqnarray}

It was shown in \cite[Corollary 4.6]{HB2} that an $\O$-operator $T:M\rightarrow \L$ on a Lie conformal algebra $\L$ associated with a representation $(M,\rho)$ also induces a left-symmetric conformal algebra structure on $M$ by
\begin{eqnarray}
m\circ_\lambda n=\rho (T(m))_\lambda n, ~~ \forall ~~m,n\in M.
\end{eqnarray}
If $\L=\A^L$ is the commutator Lie conformal algebra of an associative conformal algebra $\A$ and the representation of $\L$ on $M$ is induced from the conformal $\A$-bimodule structure on $M$ (see \eqref{L2}), then the two left-symmetric conformal algebra structures  above are exactly the same.

It was also shown in \cite[Proposition 2.2]{HB2} that if $(\A,\circ_\la)$ is a left-symmetric conformal algebra, then the $\la$-bracket
\begin{eqnarray}\label{3-00}
[a_\la b]:=a\circ_\la b-b\circ_{-\la-\partial}a, \ \forall\ a,b\in\A
\end{eqnarray}
makes $\A$ into a Lie conformal algebra, denoted by $\mathfrak{g}(\A)$. Now let $(E,\succ_\lambda,\prec_\lambda)$ be a dendriform conformal algebra. It follows from Proposition \ref{P10} that there exists a left-symmetric conformal algebra structure $(E,\circ_\la)$, which will further induce a Lie conformal algebra structure $\mathfrak{g}(E)$ by \eqref{3-00}. On the other hand, it follows from Proposition \ref{p2} that $(E,\succ_\lambda,\prec_\lambda)$ induces an associative conformal algebra $E_{ass}$. It is easy to see that the commutator Lie conformal algebra $E_{ass}^L$ of $E_{ass}$ is the same as $\mathfrak{g}(E)$.

At the end of this section, we introduce compatible $\O$-operators. Let $T_1$ and $T_2$ be two $\O$-operators on an associative conformal algebra $\A$ with respect to a conformal $\A$-bimodule $M$. They are said to be {\it compatible} if $T_1+T_2$ is again an $\O$-operator. It is not difficult to see that this requirement is equivalent to the following condition
\begin{align}\label{3-7}
T_1(m)_\la T_2(n)+T_2(m)_\la T_1(n)=T_1(T_2(m)_\la n+m_\la T_2(n))+T_2(T_1(m)_\la n+m_\la T_1(n)),
\end{align}
for all $m,n\in M.$

Note that \eqref{3-7} depends linearly on $T_1$ and $T_2$. Hence, if $T_1,\cdots, T_k$ are
$\O$-operators of $\A$ with respect to $M$ and $T_1$ is compatible with $T_2,\cdots, T_k$, then it is compatible with any linear combination of them. If $T_1,\cdots, T_k$ are pairwise compatible, then any two linear combinations of them are compatible.

Let $(E,\succ^1_\lambda,\prec^1_\lambda)$ and $(E,\succ^2_\lambda,\prec^2_\lambda)$ be two dendriform conformal algebra structures on the vector space $E$. They are said to be {\it compatible} if $(E,\succ^1_\lambda+\succ^2_\lambda,\prec^1_\lambda+\prec^2_\lambda)$  still forms a dendriform conformal algebra. Dendriform conformal algebra structures induced from
compatible $\O$-operators are compatible, as the following theorem shows.
\begin{theo}\label{th3} Let $T_1, T_2:M\rightarrow \A$ be two $\O$-operators on an associative conformal algebra $\A$ with respect to a conformal $\A$-bimodule $M$. If $T_1$ and $T_2$
are compatible, then $(M, \succ^{T_1}_\la , \prec^{T_1}_\la)$ and $(M, \succ^{T_2}_\la , \prec^{T_2}_\la)$ are compatible dendriform conformal algebras, where
\begin{align*}
m\succ^{T_1}_\la  n=T_1(m)_\lambda n, ~~ m\prec^{T_1}_\la n=m_\lambda T_1(n),~~  m\succ^{T_2}_\la  n=T_2(m)_\lambda n, ~~ m\prec^{T_2}_\la n=m_\lambda T_2(n),
\end{align*}
for all $m,n\in M.$ If, in addition, both $T_1$ and $T_2$ are invertible, then $(\A, \succ^1_\lambda,\prec^1_\lambda)$ and
$(\A, \succ^2_\lambda,\prec^2_\lambda)$ are compatible dendriform conformal algebras, where
\begin{align*}
a\succ^1_\lambda b= T_1(a_\la T_1^{-1}(b)),~~  a \prec^1_\lambda b = T_1(T_1^{-1}(a)_\la b), ~~  a\succ^2_\lambda b= T_2(a_\la T_2^{-1}(b)),~~  a \prec^2_\lambda b = T_2(T_2^{-1}(a)_\la b),
\end{align*}
for all $a,b\in \A$.
\end{theo}
\begin{proof} By Proposition \ref{p3}, it suffices to prove the compatibility. This directly follows from \eqref{3-7}
 and the definition of compatibility of dendriform conformal algebras.
\end{proof}

\section{Twisted Rota-Baxter operators and conformal NS-algebras}

In this section, we construct a twisted version of Section 3. First, we introduce the notion of a twisted Rota-Baxter operator, which is a generalization of $\O$-operators and characterized by a Hochschild $2$-cocycle. Second, we construct a new algebraic structure that is related to twisted Rota-Baxter operators in the same way that dendriform conformal algebras are related to $\O$-operators. We call such algebras as conformal NS-algebras.

\begin{defi}{\rm Let $M$ be a conformal bimodule over an associative conformal algebra $\A$, $T:M\rightarrow\A$ a $\C[\partial]$-module homomorphism and $\varphi_\lambda$ a $2$-cocycle in $C^2(\A, M)$. Then $T$ is called a twisted Rota-Baxter operator or simply $\varphi$-Rota-Baxtor operator if the condition
\begin{eqnarray}\label{twist-O}
T(m)_\lambda T(n)=T\big(T(m)_\lambda n+m_\lambda T(n)+\varphi_\lambda(T(m),T(n))\big)
\end{eqnarray} is satisfied for all $m,n\in M.$
}\end{defi}

Obviously, an $\O$-operator $T:M\rightarrow\A$ is a special twisted Rota-Baxter operator in which $\varphi=0$.
Let $\varphi_\lambda$ be any $2$-cocycle in $C^2(\A, M)$. By Theorem \ref{th10}, it corresponds
a conformal null extension $\A\oplus_{\varphi} M$ of $\A$ by means of $M$, and the associative $\lambda$-multiplication on
$\A\oplus_{\varphi} M$ is defined by \eqref{twist-sum}. Similarly to Proposition \ref{p4}, we consider the graph of $T$ and obtain the following result with a similar proof:
\begin{prop}\label{p11} $T:M\rightarrow\A$ is a $\varphi$-Rota-Baxtor operator if and only if $Gr(T)$ is a subalgebra of  $\A\oplus_{\varphi} M$.
\end{prop}

From the isomorphism $Gr(T)\cong M$, we know that $T$ induces an associative $\lambda$-multiplication
on $M$. The induced $\lambda$-multiplication of $M$ has the form
\begin{align}\label{3-1-2}
m\star^\varphi_\lambda n=T(m)_\lambda n+m_\lambda T(n)+\varphi_\lambda(T(m),T(n)), ~~ \forall ~m,n\in M.
\end{align}
We denote the new associative conformal algebra $(M,\star^\varphi_\lambda)$ by $M^\varphi_{ass}$. It is obvious that $T$ is an algebra homomorphism:
\begin{align}\label{4-12}
 T(m\star^\varphi_\lambda n)=T(m)_\la T(n), \ \forall \ m,n\in M^\varphi_{ass}.
\end{align}

\begin{prop} Let $M$ be a conformal bimodule over an associative conformal algebra $\A$. For any $2$-cocycle $\varphi\in C^2(\A,M)$ and $1$-cochain $h\in C^1(\A,M)$, we have an isomorphism of associative conformal algebras:
\begin{align*}
\A\oplus_{\varphi} M\cong \A\oplus_{\varphi+\mathbf{d} h} M.
\end{align*}
\end{prop}
\begin{proof} Define a $\C[\partial]$-module homomorphism $\psi_h:\A\oplus_{\varphi} M\rightarrow \A\oplus_{\varphi+\mathbf{d} h} M$ by
\begin{align}\label{4-2}
\psi_h(a,m)=(a,m-h(a)), \ \ \forall \ \, (a,m)\in \A\oplus_{\varphi} M.
\end{align}
Then we have
\begin{align*}
\psi_h\big((a,m)\circ^\varphi_\lambda(b,n)\big)&\stackrel{\eqref{twist-sum}}{=}\psi_h\big(a_\lambda b, a_\lambda n+m_\lambda b+\varphi_\lambda(a,b)\big)\\
&\stackrel{\eqref{4-2}}{=}\big(a_\lambda b, a_\lambda n+m_\lambda b+\varphi_\lambda(a,b)-h(a_\lambda b)\big)\\
&\stackrel{\eqref{1-cocycle}}{=}\big(a_\lambda b, a_\lambda n+m_\lambda b+\varphi_\lambda(a,b)+(\mathbf{d}h)_\la(a,b)-a_\lambda h(b)-h(a)_\la b\big)
\\&\stackrel{\eqref{twist-sum}}{=}(a,m-h(a))\circ^{\varphi+\mathbf{d}h}_\la(b,n-h(b))\\
&\stackrel{\eqref{4-2}}{=}\psi_h(a,m)\circ^{\varphi+\mathbf{d}h}_\la \psi_h(b,n).
\end{align*}
The fact that $\psi_h$ is invertible follows by exhibiting the explicit inverse $\psi^{-1}_h(a,m)=(a,m+h(a))$, for all $(a,m)\in \A\oplus_{\varphi+\mathbf{d} h} M.$
This ends the proof.
\end{proof}

\begin{ex} {\rm Let $\omega:\A\rightarrow M$ be an invertible 1-cochain in $C^1(\A, M)$. Then the inverse $\omega^{-1}$ is a twisted Rota-Baxter operator, and in this case, $\varphi=-d\omega$. In fact, putting $T=\omega^{-1}$, the condition \eqref{twist-O} is equivalent to
\begin{eqnarray*}
\omega(T(m)_\lambda T(n))=T(m)_\lambda n+m_\lambda T(n)+\varphi_\lambda(T(m),T(n)), ~~ \forall ~~ m,n\in M.
\end{eqnarray*}
This is the same as
\begin{align*}
\varphi_\lambda(T(m),T(n))&\ \, =-T(m)_\lambda n-m_\lambda T(n)+\omega(T(m)_\lambda T(n)),\\
&\stackrel{\eqref{1-cocycle}}{=} -(\mathbf{d}\omega)_\lambda(T(m),T(n)),
\end{align*}
for all $m,n\in M.$
%Recall that the coboundary of $\omega$ is defined by $(d\omega)_\lambda(a,b)=a_\lambda\omega(b)-\omega(a_\lambda b)+\omega(a)_\lambda(b)$ for all $a,b\in\A$ (see \cite{KK}).
%Hence, $\varphi_\lambda(T(m),T(n))=(d\omega)_\lambda(T(m),T(n))$.
}\end{ex}
\begin{ex} {\rm
Let $\varphi_\la\in C^2(\A,\A)$ be defined by
\begin{align}
\varphi_\lambda(a,b)=-a_\lambda b, ~~ \forall~~a,b\in\A.
\end{align}
Obviously, $\varphi$ is a 2-cocycle. Then the identity map ${\rm id}: \A\rightarrow \A$ is a $\varphi$-Rota-Baxter operator.}\end{ex}

 %Now we introduce the following notion.

\begin{defi} {\rm Let $\A$ be an associative conformal algebra. A $\C[\partial]$-module homomorphism  $R:\A\rightarrow \A$ is called a Reynolds operator of $\A$ if the condition
\begin{eqnarray}\label{R}
R(a)_\lambda R(b)=R\big(R(a)_\lambda b+a_\lambda R(b)-R(a)_\lambda R(b)\big)
\end{eqnarray} is satisfied for all $a,b\in \A.$
}\end{defi}

Notice that the last term  $-R(a)_\lambda R(b)$ in \eqref{R} is the associative $\lambda$-multiplication on $\A$, which is a 2-cocycle. Therefore each Reynolds operator $R$ can be seen as a twisted Rota-Baxter operator. It follows from \eqref{3-1-2} that $R$ induces a new associative conformal algebra structure on $\A$ by
\begin{align}
a\circ_\la^R b=R(a)_\lambda b+a_\lambda R(b)-R(a)_\lambda R(b),
\end{align}
for all $a,b\in\A.$ We denote this associative conformal algebra by $\A^R$. By \eqref{R}, $R$ is an algebra homomorphism from $\A^R$ to $\A$. Further, if $R$ is invertible, then it follows from \eqref{R} that
\begin{eqnarray}\label{R-1}
R^{-1}(a_\lambda b)=R^{-1}(a)_\lambda b+a_\lambda R^{-1}(b)-a_\lambda b,
\end{eqnarray}
for all $a,b\in\A$. This implies that $(R^{-1}-{\rm id})(a_\lambda b)=(R^{-1}-{\rm id})(a)_\lambda b+a_\lambda (R^{-1}-{\rm id})(b)$. Hence $R^{-1}-{\rm id}:\A\rightarrow \A$ is a derivation. Conversely, if $d:\A\rightarrow\A$ is a derivation such that ${\rm id}+d$ is invertible, then $({\rm id}+d)^{-1}$ is a Reynolds operator of $\A$. Even if ${\rm id}+d$
is not invertible but the infinite sum $({\rm id}+d)^{-1}=\sum_{n=0}^\infty (-1)^nd^n$ converges pointwise, then $({\rm id}+d)^{-1}$ is a Reynolds operator of $\A$. A more precise
statement is given below by a verbatim repetition of the proof of \cite[Proposition 2.8]{ZGG} in terms of $\la$-multiplication.

\begin{prop}
Let $\A$ be an associative conformal algebra with a derivation $d$. If the series $\sum_{n=0}^\infty (-1)^nd^n (x)$ is convergent for all $x\in\A$, then $R:=\sum_{n=0}^\infty (-1)^nd^n$ is a Reynolds operator of $\A$.
\end{prop}
It follows from the above proposition that if $d$ is a nilpotent derivation (more generally, a locally nilpotent derivation) on $\A$, then $R=\sum_{n=0}^\infty (-1)^nd^n$ is a Reynolds operator of $\A$.

We see from Proposition \ref{p3} that $\mathcal{O}$-operators induce dendriform conformal algebra
structures. In the following, we will show a similar result with respect to twisted Rota-Baxter operators.
We need the following concept.

\begin{defi}{\rm Let $\mathcal{N}$ be a $\C[\partial]$-module equipped with three binary $\lambda$-multiplications
$\succ_\lambda, \prec_\lambda$ and $\vee_\lambda$. Then $\mathcal{N}$ is called a conformal NS-algebra, if $\succ_\lambda, \prec_\lambda$ and $\vee_\lambda$ are conformal sesquilinear maps, and satisfy the following axioms for all $x,y,z\in\mathcal{N}$:
\begin{align}
x\succ_\lambda(y\succ_\mu z)&=(x\times_\lambda y)\succ_{\lambda+\mu}z,\label{NS1}\\
x\prec_\lambda (y\times_\mu z)&=(x\prec_\lambda y)\prec_{\lambda+\mu}z,\label{NS2}\\
x\succ_\lambda (y\prec_\mu z)&=(x\succ_\lambda y)\prec_{\lambda+\mu}z,\label{NS3}\\
x\succ_\lambda (y\vee_\mu z)- (x\times_\lambda y)\vee_{\lambda+\mu}z&=(x\vee_\lambda y)\prec_{\lambda+\mu}z-x\vee_\lambda (y\times_\mu z),\label{NS4}
\end{align}
where $\times_\lambda$ is defined as
\begin{align}\label{NS5}
x\times_\lambda y=x\succ_\lambda y+x\prec_\lambda y+x\vee_\lambda y.
\end{align}
}\end{defi}

The basic property of usual dendriform conformal algebras is satisfied on conformal
NS-algebras.

\begin{prop}\label{p5}  Let $\mathcal{N}$ be a conformal
NS-algebra. Then $(\mathcal{N}, \times_\lambda)$ forms an associative conformal algebra, where $\times_\lambda$ is defined by \eqref{NS5}.
\end{prop}
\begin{proof}
Straightforward.
\end{proof}

The following theorem reveals a close relation between twisted Rota-Baxter operators and conformal NS-algebras.

\begin{theo}\label{th7} Assume that $M$ is a conformal bimodule over an associative conformal algebra $\A$, $\varphi_\lambda$ is a $2$-cocycle in $C^2(\A, M)$, and $T:M\rightarrow \A$ is a
$\varphi$-Rota-Baxter operator. Then $M$ becomes a conformal NS-algebra under the following three $\la$-multiplications:
\begin{eqnarray}\label{6-6}
m\succ_\la n=T(m)_\la n, ~~ m\prec_\la n=m_\la T(n), ~~ m\vee_\la n=\varphi_\la(T(m),T(n)),
\end{eqnarray}
where $m,n\in M.$
\end{theo}
\begin{proof}
 It is easy to see that $\succ_\lambda, \prec_\lambda$ and $\vee_\lambda$ are conformal  sesquilinear maps. Relation \eqref{NS3} is easy to check. For any $m,n,l\in M$, we have
\begin{align*}
m\prec_\la(n\times_\mu l)&=m\prec_\la\big(T(n)_\mu l+n_\mu T(l)+ \varphi_\mu(T(n),T(l))\big)\\
&=m_\la T\big(T(n)_\mu l+n_\mu T(l)+ \varphi_\mu(T(n),T(l))\big)\\
&=m_\la(T(n)_\mu T(l))=(m_\la T(n))_{\la+\mu}T(l)\\
&=(m\prec_\la n)\prec_{\la+\mu} l.
\end{align*}
This proves \eqref{NS2}. Relation \eqref{NS1} can be similarly obtained. It is left to show \eqref{NS4}. Because $\varphi$ is a Hochschild 2-cocycle, we have the following cocycle condition:
\begin{eqnarray*}
0=T(m)_\lambda \varphi_\mu (T(n),T(l))-\varphi_{\lambda+\mu}(T(m)_\lambda T(n),T(l))+\varphi_{\lambda}(T(m), T(n)_\mu T(l))-\varphi_{\lambda}(T(m),T(n))_{\lambda+\mu}T(l),
\end{eqnarray*}
for all $m,n,l\in M.$ This, together with \eqref{6-6}, gives
\begin{eqnarray}\label{7-7}
0=m\succ_\lambda (n\vee_\mu l)-\varphi_{\lambda+\mu}(T(m)_\lambda T(n),T(l))+\varphi_{\lambda}(T(m), T(n)_\mu T(l))-(m\vee_{\lambda}n)\prec_{\lambda+\mu} l.
\end{eqnarray}
On the other hand, we have
\begin{eqnarray*}
T(m)_\lambda T(n)=T\big( T(m)_\la n+m_\la T(n)+\varphi_\la(T(m),T(n)) \big)=T(m\succ_\la n+m\prec_\la n+m\vee_\la n)=T(m\times_\la n).
\end{eqnarray*}
It follows that
\begin{align*}
 \varphi_{\lambda+\mu}(T(m)_\lambda T(n),T(l))&=\varphi_{\lambda+\mu}(T(m\times_\la n),T(l))=(m\times_\la n)\vee_{\la+\mu}l,\\
 \varphi_{\lambda}(T(m), T(n)_\mu T(l))&=\varphi_{\lambda}(T(m), T(n\times_\mu l))=m\vee_\la(n\times_\mu l).
\end{align*}
Plugging this back into \eqref{7-7}, we obtain \eqref{NS4}. This ends the proof.
\end{proof}
\begin{remark}{\rm \label{coro2}Proposition \ref{p5} and Theorem \ref{th7} also imply that $(M,\star^\varphi_\la)$ is an associative conformal algebra, where the associative $\la$-multiplication $\star^\varphi_\la$ is defined by \eqref{3-1-2}. And $T$ is an associative conformal algebra homomorphism from $(M,\star^\varphi_\la)$ to $\A$.}
\end{remark}

\begin{prop}\label{p18} If $T:M\rightarrow \A$ is a
$\varphi$-Rota-Baxter operator, then $\A$ becomes a conformal $M^\varphi_{ass}$-bimodule
by the following $\lambda$-actions:
\begin{align}\label{4-11}
m\cdot^\varphi_\lambda a=T(m)_\lambda a-T\big(m_\lambda a+\varphi_\la(T(m),a)\big), ~~ a\cdot^\varphi_\lambda m=a_\lambda T(m)-T\big(a_\lambda m+\varphi_\la(a,T(m))\big),
\end{align}
where $m\in M^\varphi_{ass}$ and $a\in\A.$
\end{prop}
\begin{proof} It is easy to see that the two $\lambda$-actions defined by \eqref{4-11} are conformal sesquilinear maps. For any $m,n\in M$ and $a,b\in\A,$ we have
\begin{align*}
m\cdot^\varphi_\lambda (n\cdot^\varphi_\mu a)=&\, T(m)_\la(n\cdot^\varphi_\mu a)-T\big(m_\la(n\cdot^\varphi_\mu a)+\varphi_\la(T(m),n\cdot^\varphi_\mu a)\big)\\
=&\,T(m)_\la\big(T(n)_\mu a-T(n_\mu a+\varphi_\mu(T(n),a))\big)
-T\big(m_\la (T(n)_\mu a-T(n_\mu a+\varphi_\mu(T(n),a))\big)\\
&-T\varphi_\la\big(T(m), T(n)_\mu a-T(n_\mu a+\varphi_\mu(T(n),a))\big)\\
=&\,T(m)_\lambda(T(n)_\mu a)-T\big( T(m)_\la(n_\mu a)+\cancel{m_\la T(n_\mu a)}+\cancel{\varphi_\la(T(m),T(n_\mu a))} \big)\\
&-T\big(T(m)_\la\varphi_\mu(T(n),a)+\cancel{m_\la T\varphi_\mu(T(n),a)}+\cancel{\varphi_\la(T(m),T\varphi_\mu(T(n),a))}\big)\\
&-T\big( m_\la(T(n)_\mu a)-\cancel{m_\la T(n_\mu a)}-\cancel{m_\la T\varphi_\mu(T(n),a)}\big)\\
&+T\varphi_\la\big(T(m),-T(n)_\mu a+\cancel{T(n_\mu a)}+\cancel{T\varphi_\mu(T(n),a)}\big)\\
=&\,T(m)_\lambda(T(n)_\mu a)-T\big(m_\lambda (T(n)_\mu a)+T(m)_\lambda (n_\mu a)\big)\\
&-T\big(T(m)_\la\varphi_\mu(T(n),a)+\varphi_\la(T(m),T(n)_\mu a)\big).
\end{align*}
On the other hand, by \eqref{3-1-2} and \eqref{4-12}, we have
\begin{align*}
(m\star^\varphi_\lambda n)\cdot^\varphi_{\lambda+\mu} a=&T(m\star^\varphi_\lambda n)_{\lambda+\mu} a-T\big((m\star^\varphi_\lambda n)_{\lambda+\mu} a+\varphi_{\la+\mu}(T(m\star^\varphi_\lambda n),a)\big)\\
=&(T(m)_\lambda T(n))_{\lambda+\mu} a-T\big((T(m)_\lambda n+m_\lambda T(n)+\varphi_\la(T(m),T(n)))_{\lambda+\mu} a+\varphi_{\la+\mu}(T(m)_\la T(n),a)\big).
%&-T\big(\varphi_\la(T(m),T(n))_{\la+\mu}a+\varphi_{\la+\mu}(T(m)_\la T(n),a)\big).
\end{align*}
Then, by the fact that $M$ is a conformal $\A$-bimodule and $\varphi$ is a 2-cocycle in $C^2(\A,M)$, we obtain $$m\cdot^\varphi_\lambda (n\cdot^\varphi_\mu a)=(m\star^\varphi_\lambda n)\cdot^\varphi_{\lambda+\mu} a.$$
Similarly, we can obtain $(a\cdot^\varphi_\lambda m)\cdot^\varphi_{\lambda+\mu} n=a\cdot^\varphi_\lambda (m\star^\varphi_{\mu}n)$ and $m\cdot^\varphi_\lambda (a\cdot^\varphi_\mu n)=(m\cdot^\varphi_\lambda a)\cdot^\varphi_{\lambda+\mu} n$. Hence $\A$ is a conformal $M^\varphi_{ass}$-bimodule.
\end{proof}

Let $\mathcal{N}$ be a conformal NS-algebra. We denote the associated associative conformal algebra $(\mathcal{N}, \times_\la)$ by $\mathcal{N}_{ass}$. A conformal $\mathcal{N}_{ass}$-bimodule
structure on $\mathcal{N}$ is well-defined by
\begin{align}
x_\la t= x \succ_\la t,\ \ t_\la x= t\prec _\la x,
\end{align}
where $x\in \mathcal{N}_{ass}$ and $t\in\mathcal{N}$.

\begin{prop}\label{p12} Under the assumptions above, for any $x,y\in \mathcal{N}_{ass}$, define $\varphi_\la(x,y)=x\vee_\la y$. Then $\varphi_\la$ is a 2-cocycle in $C^2(\mathcal{N}_{ass},\mathcal{N})$ and thus the identity map
${\rm {id}}: \mathcal{N}\rightarrow \mathcal{N}_{ass}$ is a $\varphi$-Rota-Baxter operator.
\end{prop}
\begin{proof} The cocycle condition of $\varphi$
 is the same as \eqref{NS4}. %Hence the first assertion is valid. The second assertion is obvious.
 \end{proof}
 %In fact,\begin{align*}
%&x_\lambda \varphi_\mu (y,z)-\varphi_{\lambda+\mu}(x\times_\lambda y,z)+\varphi_{\lambda}(x, y\times_\mu z)-\varphi_{\lambda}(x,y)_{\lambda+\mu}z\\=&x\succ_\la(y\vee_\mu z )-(x\times_\la y)\vee_{\la+\mu} z+x\vee_\la(y\times_\mu z )-(x\vee_\la y)\prec_{\la+\mu}z\\=&0.
%\end{align*}

\begin{remark} {\rm Under the setting in Proposition \ref{p12}, it follows from Theorem \ref{th7} that there is a new conformal NS-algebra structure on $\mathcal{N}$, which coincides exactly with the original one.}
\end{remark}

Let $M$ be a conformal bimodule over an associative conformal algebra $\A$. A 2-cocycle $\varphi_\la$ in $C^2(\A,M)$ is said to be {\it commutative} if the condition
\begin{align}
\varphi_\la(a,b)=\varphi_{-\la-\partial}(b,a)
\end{align}
is satisfied for all $a,b\in\A.$

We have the following result with a similar proof of Theorem \ref{th2}:
\begin{theo} \label{p13} Suppose that $M$ is a conformal bimodule over an associative conformal algebra $\A$ and $\varphi_\la$ is a commutative 2-cocycle in $C^2(\A,M)$. If $T:M\rightarrow \A$ is a $\varphi$-Rota-Baxter operator, then $T$ is also an $\mathcal{O}$-operator on the commutator Lie conformal algebra $\A^L$ with respect to the representation $(M,\rho)$.
\end{theo}

Let $T:M\rightarrow \A$ be a
$\varphi$-Rota-Baxter operator. It follows from \eqref{3-1-2} (see also Remark \ref{coro2}) that $(M,\star^\varphi_\la)$ becomes an associative conformal algebra. If, in addition, $\varphi_\la$ is commutative, then the commutator Lie conformal algebra structure associated with
$(M,\star^\varphi_\la)$ is the same as the untwisted one (cf. \eqref{untwist}).

\section{Nijenhuis operators of associative conformal algebras}
In this section, we introduce a conformal analog of associative Nijenhuis operators, and enumerate main properties. Further, we present connections between deformations and Nijenhuis operators of associative conformal algebras.

\subsection{Definition and properties of Nijenhuis operators}

We first introduce the notion of a Nijenhuis operator for arbitrary associative conformal algebras.
\begin{defi} {\rm Let $\A$ be an associative conformal algebra. A $\C[\partial]$-module homomorphism $N:\A\rightarrow \A$ is called a Nijenhuis operator of $\A$ if the condition
\begin{eqnarray}\label{N}
N(a)_\lambda N(b)=N\big(N(a)_\lambda b+a_\lambda N(b)-N(a_\lambda b)\big)
\end{eqnarray} is satisfied for all $a,b\in \A.$
}\end{defi}

Obviously, the identity map $\rm {id}$ is  a Nijenhuis operator of $\A$.

\begin{defi}{\rm (see \cite{HB2}) Let $\A$ be an associative conformal algebra and $T:\A\rightarrow \A$ a $\C[\partial]$-module homomorphism. For $q\in\C$, if there holds
\begin{eqnarray}
T(a)_\lambda T(b)=T\big(T(a)_\lambda b+a_\lambda T(b)+q a_\la b\big),
\end{eqnarray}
for all $a,b\in\A$, then $T$ is called a Rota-Baxter operator of weight $q$ on $\A$.
}\end{defi}
The following proposition describes close interrelations between Nijenhuis operators and Rota-Baxter
operators, and the proof is straightforward.
\begin{prop} \label{p19}Let $N:\A\rightarrow \A$ be a $\C[\partial]$-module homomorphism over an associative conformal algebra $\A$.
\begin{itemize}
\item[(i)] If $N^2 =0$, then $N$ is a Nijenhuis operator if and only if $N$ is a Rota-Baxter
operator of weight 0.
\item[(ii)] If $N^2 =N$, then $N$ is a Nijenhuis operator if and only if $N$ is a Rota-Baxter
operator of weight $-1$.
\item[(iii)] If $N^2 ={\rm id}$, then $N$ is a Nijenhuis operator if and only if $N\pm{\rm id}$ is a Rota-Baxter
operator of weight $\mp 2$.
\end{itemize}
\end{prop}

 \begin{ex}{\rm Let $T:M\rightarrow \A$ be an  $\mathcal{O}$-operator on an associative conformal algebra $\A$ with respect to a conformal $\A$-bimodule $M$. By Proposition \ref{p3-3}, the lift $\hat{T}$ is a Rota-Baxter operator of weight 0 on  $\A\oplus_0 M$. Obviously, $\hat{T}^2=0$. Hence, by Proposition \ref{p19} (i),  $\hat{T}$ is a Nijenhuis operator of $\A\oplus_0 M$.}
\end{ex}

The following theorem says that Nijenhuis operators on associative conformal algebras give rise to conformal NS-algebra structures.

\begin{theo}\label{p6}  Let $N$ be a Nijenhuis operator over an associative conformal algebra $\A$.  For all $a,b\in\A$,
define three $\lambda$-multiplications on $\A$ by
\begin{align}\label{5-5}
a\succ_\lambda b=N(a)_\lambda b,~~ a\prec_\lambda b=a_\lambda N(b),~~ a\vee_\lambda b=-N(a_\lambda b).
\end{align}
Then $(\A,\succ_\lambda, \prec_\lambda,\vee_\lambda)$ is a conformal NS-algebra.
\end{theo}
\begin{proof} Obviously, $\succ_\lambda, \prec_\lambda$ and $\vee_\lambda$ are conformal sesquilinear maps. It follows from \eqref{N} and \eqref{5-5} that
\begin{align}
N(a)_\lambda N(b)=N(a\times_\lambda b), \ \forall \ a,b\in\A,
\end{align}
where $a\times_\lambda b=a\succ_\lambda b+a\prec_\lambda b+a\vee_\lambda b$.
Hence, for any $ a,b,c\in\A$, we have
\begin{align*}
a\succ_\lambda(b\succ_\mu c)=N(a)_\lambda(N(b)_\mu c)=(N(a)_\lambda N(b))_{\lambda+\mu} c=N(a\times_\lambda b)_{\lambda+\mu} c=(a\times_\lambda b)\succ_{\lambda+\mu} c.
\end{align*}
This proves \eqref{NS1}. Relations \eqref{NS2} and \eqref{NS3} can be similarly proved. To prove \eqref{NS4}, we compute, respectively,
\begin{align}
a\succ_\lambda(b& \vee_\mu c)- (a\times_\lambda b)\vee_{\lambda+\mu}c\nonumber\\
=&-N(a)_\la N(b_\mu c)+N\big((a\times_\lambda b)_{\lambda+\mu}c\big)\nonumber\\
=&-N\big(N(a)_\la(b_\mu c)+a_\la N(b_\mu c)-N(a_\la(b_\mu c))\big)
 +N\big((N(a)_\la b+a_\la N(b)-N(a_\la b))_{\lambda+\mu}c\big)\nonumber\\
%=&-N\big(a_\la N(b_\mu c)-N(a_\la(b_\mu c))\big)
 %+N\big(a_\la N(b)-N(a_\la b))_{\lambda+\mu}c\big)\nonumber\\
=&-N(a_\la N(b_\mu c))+N^2(a_\la(b_\mu c))+N(
 a_\la (N(b)_\mu c))-N\big(N(a_\la b)_{\lambda+\mu}c\big),\label{NS6}
\end{align}
and
\begin{align}
(a\vee_\lambda b)&\prec_{\lambda+\mu}c-a\vee_\lambda (b\times_\mu c)\nonumber\\
&=-N(a_\la b)_{\la+\mu}N(c)+N\big(a_\la(b\times_\mu c)\big)\nonumber\\
&=-N\big(N(a_\la b)_{\la+\mu}c+ (a_\la b)_{\la+\mu}N(c)-N((a_\la b)_{\la+\mu}c)\big)+N\big(a_\la(N(b)_\mu c+b_\mu N(c)-N(b_\mu c))\big)\nonumber\\
&=-N\big(N(a_\la b)_{\la+\mu}c\big)+N^2((a_\la b)_{\la+\mu}c)+N\big(a_\la(N(b)_\mu c-N(b_\mu c))\big)\label{NS7}.
%=&-N\big(N(a_\la b)_{\la+\mu}c-N((a_\la b)_{\la+\mu}c)-a_\la(N(b)_\mu c-N(b_\mu c))\big)\nonumber\\
%=&-N\big(N(a_\la b)_{\la+\mu}c-N(a_\la(b_\mu c))-(a_\la N(b))_{\la+\mu} c+a_\la N(b_\mu c)\big).
\end{align}
Comparing \eqref{NS6} with \eqref{NS7} gives \eqref{NS4}. This completes the proof.
 %$a\succ_\lambda (b\vee_\mu c)- (a\times_\lambda b)\vee_{\lambda+\mu}c=(a\vee_\lambda b)\prec_{\lambda+\mu}c-a\vee_\lambda (b\times_\mu c)$.
\end{proof}

Combining Proposition \ref{p5} and Theorem \ref{p6}, we have the following corollary.

\begin{coro}\label{coro1} Let $N$ be a Nijenhuis operator over an associative conformal algebra $\A$.  Define
\begin{eqnarray}\label{N2}
a \circ^N_\lambda b:=N(a)_\lambda b+a_\lambda N(b)-N(a_\lambda b), \ \forall \  a,b\in \A.
\end{eqnarray}
Then $(\A, \circ^N_\lambda)$ forms a new associative conformal algebra, denoted by $\A^N$. Further, $N$ is an algebra homomorphism from $\A^N$ to the original associative conformal algebra $\A$:
\begin{eqnarray}\label{N3}
N(a \circ^N_\lambda b)=N(a)_\la N(b), \ \forall \ a,b\in\A.
\end{eqnarray}
\end{coro}

Assume that $N:\A\rightarrow \A$ is a $\C[\partial]$-module homomorphism over an associative conformal algebra $\A$. In the following we denote the $\la$-product on $\A$ by $\theta_\la$, i.e., $\theta_\la(a,b)=a_\la b$ for all $a,b\in\A.$ The map
\begin{eqnarray}\label{ast}
\theta^N_\la:(a,b)\mapsto a \circ^N_\lambda b=N(a)_\lambda b+a_\lambda N(b)-N(a_\lambda b), \ \forall \  a,b\in \A
\end{eqnarray}
is conformal sesquilinear and therefore it defines a new algebra structure on $\A$.
Then we obtain a 2-cochain $\varphi^N_\la$ in $C^2(\A,\A)$ of the form
\begin{eqnarray}\label{ast2}
\varphi^N_\lambda(a,b)=N(a)_\lambda N(b)-N(a\circ^N_\lambda b), \ \forall \  a,b\in \A.
\end{eqnarray}
It is obvious that $\varphi^N_\lambda=0$ if and only if $N$ is a Nijenhuis operator of $\A$.% Moreover, if $N$ is a Nijenhuis operator of $\A$, then it follows Propositions \ref{p5} and \ref{p6} that $\A^N:=(\A,\theta^N_\la)$ forms a new associative conformal algebra and $N$ is an algebra homomorphism from $\A^N$  to $\A$.

\begin{theo}\label{th1} Under the assumptions above, $(\A,\theta^N_\lambda)$ forms an associative conformal algebra if and only if $\varphi^N_\la$ is a 2-cocycle in $C^2(\A,\A)$, i.e.,
\begin{align}
({\mathbf{d}} \varphi^N)_{\la,\mu}(a,b,c):=a_\lambda \varphi^N_\mu (b,c)-\varphi^N_{\lambda+\mu}(a_\lambda b,c)+\varphi^N_{\lambda}(a, b_\mu c)-\varphi^N_{\lambda}(a,b)_{\lambda+\mu}c=0,
\end{align}
for all $a,b,c\in\A$. If this is the case, then  $\theta^N_\la$ is an associative $\la$-product compatible with $\theta_\la$, i.e., the maps $\theta_\la+q\theta^N_\la$ are associative for all $q\in\C$.
%In particular, if $N$ a Nijenhuis operator of $\A$, then $\theta^N_\la$ is an associative $\la$-product on $\A$  compatible with $\theta_\la$.
\end{theo}
\begin{proof}First we show there holds
\begin{align}\label{2-2-2}
(a\circ^N_\la b)_{\la+\mu}c+(a_\la b)\circ^N_{\la+\mu}c=a_\la(b\circ^N_\mu c)+a\circ^N_\la(b_\mu c)
\end{align}
for all $a,b,c\in\A.$  This is straightforward:
\begin{align*}
(a\circ^N_\la b)&_{\la+\mu}c+(a_\la b)\circ^N_{\la+\mu}c\\&\stackrel{\eqref{N2}}{=}\big(N(a)_\la b+a_\la N(b)-N(a_\la b)\big)_{\la+\mu}c+N(a_\la b)_{\la+\mu}c+(a_\la b)_{\la+\mu}N(c)-N((a_\la b)_{\la+\mu}c)\\
&=\big(N(a)_\la b+a_\la N(b)\big)_{\la+\mu}c+(a_\la b)_{\la+\mu}N(c)-N((a_\la b)_{\la+\mu}c)\\
&\stackrel{\eqref{ASS}}{=}a_\la (b_{\mu}N(c))+a_\la (N(b)_\mu c)+N(a)_\la (b_\mu c)-N(a_\la (b_\mu c))\\
&=a_\la (b_{\mu}N(c))+a_\la (N(b)_\mu c)-a_\la N(b_\mu c)+a_\la N(b_\mu c)+N(a)_\la (b_\mu c)-N(a_\la (b_\mu c))\\
&\stackrel{\eqref{N2}}{=}a_\la(b\circ^N_\mu c)+a\circ^N_\la(b_\mu c).
\end{align*}

Next we compute, separately,
\begin{align}\label{2-2-4}
(a\circ^N_\la b)\circ^N_{\la+\mu}c-&a\circ^N_\la (b\circ^N_{\mu}c)\nonumber\\
=&N(a\circ^N_\la b)_{\la+\mu}c+(a\circ^N_\la b)_{\la+\mu}N(c)-N((a\circ^N_\la b)_{\la+\mu}c)\nonumber\\&-N(a)_\la(b\circ^N_\mu c)-a_\la N(b\circ^N_\mu c)+N(a_\la(b\circ^N_\mu c))\nonumber\\
=&N(a\circ^N_\la b)_{\la+\mu}c-a_\la N(b\circ^N_\mu c)-N\big((a\circ^N_\la b)_{\la+\mu}c-a_\la(b\circ^N_\mu c)\big) \nonumber\\
&+\big(a_\la N(b)+N(a)_\la b-N(a_\la b)\big)_{\la+\mu}N(c)-N(a)_\la\big(N(b)_\mu c+b_\mu N(c)-N(b_\mu c)\big)\nonumber\\
=&N(a\circ^N_\la b)_{\la+\mu}c-a_\la N(b\circ^N_\mu c)-N\big((a\circ^N_\la b)_{\la+\mu}c-a_\la(b\circ^N_\mu c)\big) \nonumber\\
&+a_\la (N(b)_\mu N(c))-N(a_\la b)_{\la+\mu}N(c)-(N(a)_\la N(b))_{\la+\mu} c+N(a)_\la N(b_\mu c),
\end{align}
and
\begin{align}\label{2-2-3}
(\mathbf{d} \varphi^N)_{\la,\mu}(a,b,c)=&a_\lambda \varphi^N_\mu (b,c)-\varphi^N_{\lambda+\mu}(a_\lambda b,c)+\varphi^N_{\lambda}(a, b_\mu c)-\varphi^N_{\lambda}(a,b)_{\lambda+\mu}c\nonumber\\
=&a_\lambda(N(b)_\mu N(c)-N(b\circ^N_\la c))-(N(a_\la b)_{\la+\mu}N(c)-N((a_\la b)\circ^N_{\la+\mu}c))\nonumber\\
&+N(a)_\la N(b_\mu c)-N(a\circ^N_\la(b_\mu c))-(N(a)_{\lambda}N(b)-N(a\circ^N_\la b))_{\lambda+\mu}c\nonumber\\
=&N(a\circ^N_\la b)_{\lambda+\mu}c-a_\lambda N(b\circ^N_\la c)+a_\lambda(N(b)_\mu N(c))-N(a_\la b)_{\la+\mu}N(c)\nonumber\\
&+N(a)_\la N(b_\mu c)-(N(a)_{\lambda}N(b))_{\lambda+\mu}c+N\big((a_\la b)\circ^N_{\la+\mu}c-a\circ^N_\la(b_\mu c)\big).
\end{align}
Combining \eqref{2-2-3} with \eqref{2-2-4}, and utilizing \eqref{2-2-2}, we obtain $(a\circ^N_\la b)\circ^N_{\la+\mu}c-a\circ^N_\la (b\circ^N_{\mu}c)=(\mathbf{d} \varphi^N)_{\la,\mu}(a,b,c)$. This implies the first assertion.

As for the associativity of $\theta_\la+q\theta^N_\la$, it is exactly equivalent to \eqref{2-2-2}.
The proof is completed.
\end{proof}

\begin{remark} {\rm Relation \eqref{2-2-2} implies that the map $\theta^N_\la$ as a 2-cochain in $C^2(\A,\A)$ is exactly a 2-cocycle. Moreover, \eqref{2-2-2} holds automatically, no matter if $\theta^N_\la$ is associative or not.
Hence, if we look for a new $\la$-product $\circ_\la$ which is compatible in the sense of \eqref{2-2-2}, then this means that the new $\la$-product $\circ_\la$ is a 2-cocycle of the original associative conformal algebra. If our algebra is, for instance, the current conformal algebra ${\rm Cur}_n$ or the conformal algebra ${\rm Cend}_n$, it is proved in \cite{D2} that the second cohomology group of
${\rm Cend}_n$ and ${\rm Cur}_n$ with coefficients in any conformal bimodule is trivial, hence our  $\la$-product $\circ_\la$ has to be a coboundary, namely, of the form $\circ^N_\la$ for some $N$. This means that we have
not much freedom and, looking for compatible associative $\la$-products, we must, in
principle, work with Nijenhuis operators.}
\end{remark}

\begin{ex}{\rm Let $N:\A\rightarrow \A$ be a Nijenhuis operator on an associative conformal algebra $\A$. Then $\A$ becomes a conformal $\A^N$-bimodule by
\begin{align}
a\cdot_\la x=N(a)_\la x,\ \ x\cdot_\la a=x_\la N(a),
\end{align}
where $a\in\A^N$ and $x\in\A.$ With this bimodule, the map $\varphi_\la(a,b):=-N(a_\la b)$ is a 2-cocycle in $C^2(\A^N,\A)$. Then it is easy to see that the identity map ${\rm id}: \A\rightarrow \A^N$ is a $\varphi$-twisted Rota-Baxter operator.
}\end{ex}

\begin{lemm}\label{lemm3} Let $N:\A\rightarrow \A$ be a Nijenhuis operator on an associative conformal algebra $\A$.
For arbitrary elements $a,b\in\A$ and arbitrary nonnegative numbers $j,k\in\Z$, there holds
\begin{eqnarray}\label{3-0}
N^j(a)_\la N^k(b)-N^k(N^j(a)_\la b)-N^j(a_\la N^k(b))+N^{j+k}(a_\la b)=0.
\end{eqnarray}
If $N$ is invertible, this formula is valid for arbitrary $j,k\in\Z$.
\end{lemm}
\begin{proof} It is easy to see that \eqref{3-0} always holds for either $j=0$ or $k=0$. Now we
fix $j=1$ and prove \eqref{3-0} for arbitrary $k> 0$. For $k=1$, the formula is exactly \eqref{N}.
With the help of \eqref{N}, we get
\begin{align*}
N(a)&_\la N^{k+1}(b)-N^{k+1}(N(a)_\la b)-N(a_\la N^{k+1})+N^{k+2}(a_\la b)\\
&=N(N(a)_\la N^k(b))-N^2(a_\la N^k(b))-N^{k+1}(N(a)_\la b)+N^{k+2}(a_\la b)\\
&=N\big(N(a)_\la N^k(b)-N^{k}(N(a)_\la b)-N(a_\la N^k(b))+N^{k+1}(a_\la b)\big).
\end{align*}
By induction it follows that
\begin{align}\label{3-0-1}
N(a)_\la N^k(b)-N^{k}(N(a)_\la b)-N(a_\la N^k(b))+N^{k+1}(a_\la b)=0.
\end{align}
Now applying \eqref{3-0-1} to the element $N^j(a)$ instead of the element $a$ gives
\begin{align}\label{3-0-1-1}
N^{j+1}(a)_\la N^k(b)-N^{k}(N^{j+1}(a)_\la b)-N(N^j(a)_\la N^k(b))+N^{k+1}(N^j(a)_\la b)=0.
\end{align}
%and again relying on \eqref{N},
Then we obtain
\begin{align*}
N^{j+1}(a)_\la &N^k(b)-N^k(N^{j+1}(a)_\la b)-N^{j+1}(a_\la N^k(b))+N^{j+k+1}(a_\la b)\\
&\stackrel{\eqref{3-0-1-1}}{=}N(N^{j}(a)_\la N^k(b))-N^{k+1}(N^j(a)_\la b)-N^{j+1}(a_\la N^k(b))+N^{j+k+1}(a_\la b)\\
&\ \,=N\big( N^j(a)_\la N^k(b)-N^k(N^j(a)_\la b)-N^j(a_\la N^k(b))+N^{j+k}(a_\la b)\big).
\end{align*}
The conclusion is that the induction can be made with respect to $j$, starting from the formula \eqref{3-0-1} already proved. Thus we have proved the validity of \eqref{3-0} for arbitrary $j,k> 0$.

Suppose that $N$ is invertible. Applying $N^{-k}$ to formula \eqref{3-0} and substituting $b_1=N^k(b)$, we have
\begin{eqnarray*}
N^{-k}(N^j(a)_\la b_1)-N^j(a)_\la N^{-k}(b_1)-N^{j-k}(a_\la b_1)+N^{j}(a_\la N^{-k}(b_1))=0.
\end{eqnarray*}
As $b_1$ can be taken arbitrarily, \eqref{3-0} also holds for $k< 0$, $j> 0$. Similarly, \eqref{3-0} holds for $k> 0$, $j<0$. To prove \eqref{3-0} for both $k,j$ negative, we can apply $N^{-j-k}$ to \eqref{3-0} with putting $a_1=N^j(a)$ and $b_1=N^{k}(b)$. This ends the proof.
\end{proof}

\begin{prop}\label{p14} Let $N:\A\rightarrow \A$ be a Nijenhuis operator on an associative conformal algebra $\A$. Then for any polynomial $P(z)=\sum_{i=0}^n c_i z^i$, $P(N)$ is also a Nijenhuis operator. If $N$ is invertible, then for any $Q(z)=\sum_{i=-m}^n c_i z^i$, $Q(N)$ is also a Nijenhuis operator.
\end{prop}
\begin{proof} For arbitrary $a,b\in\A$, we have
\begin{align*}
P(N)&(a)_\la P(N)(b)-P(N)\big(P(N)(a)_\lambda b+a_\lambda P(N)(b)-P(N)(a_\lambda b)\big)\\
&=\sum_{k,j=0}^n c_j c_k\big(N^j(a)_\la N^k(b)-N^k(N^j(a)_\la b)-N^j(a_\la N^k(b))+N^{j+k}(a_\la b)\big).
\end{align*}
But the right-hand side of this equality vanishes due to \eqref{3-0}. The second statement is valid for similar reasons.
\end{proof}

\begin{lemm}\label{lemm2} Let $N:\A\rightarrow \A$ be a Nijenhuis operator on an associative conformal algebra $\A$. Then for all $a,b\in\A$ and $k,r=0,1,2,\cdots$, there holds
\begin{align}\label{3-1}
N^r(a\circ^{N^{k+r}}_\la b)=N^r(a)\circ^{N^k}_\la N^r(b),
\end{align}
i.e.,
\begin{align}\label{3-1-1}
N^r\big(N^{k+r}(a)_\la b+a_\la N^{k+r}(b)-N^{k+r}(a_\la b)\big)=N^{k+r}(a)_\la N^r(b)+N^{r}(a)_\la N^{k+r}(b)-N^k(N^{r}(a)_\la N^{r}(b)).
\end{align}
\end{lemm}
\begin{proof} The case of $r=0$ is trivial and the case of $k=0$ is equivalent to say that $N^r$ is a Nijenhuis operator of $\A$, which is valid due to Proposition \ref{p14}. Now for $r=1$ we prove
\begin{align}\label{3-2}
N(a\circ^{N^{k+1}}_\la b)=N(a)\circ^{N^k}_\la N(b)
\end{align}
holds for arbitrary $k>0$. Applying $N^k$ to \eqref{N} gives
\begin{align}\label{3-3}
N^{k+2}(a_\la b)-N^{k+1}(a_\la N(b))=N^{k+1}(N(a)_\la b)-N^{k}(N(a)_\la N(b)).
\end{align}
Using \eqref{3-3} inductively for $k:=k-1$, we end up with
\begin{align}\label{3-4}
N^{k+2}(a_\la b)-N^{k+1}(a_\la N(b))=N(N^{k+1}(a)_\la b)-N^{k+1}(a)_\la N(b).
\end{align}
In a similar way, we get
\begin{align*}
N^{k+2}(a_\la b)-N^{k+1}(N(a)_\la b)=N(a_\la N^{k+1}(b))-N(a)_\la N^{k+1}(b),
\end{align*}
which, combined with \eqref{3-3}, gives
\begin{align}\label{3-5}
N^{k+1}(a_\la N(b))-N^{k}(N(a)_\la N(b))=N(a_\la N^{k+1}(b))-N(a)_\la N^{k+1}(b).
\end{align}
Combining \eqref{3-4} and \eqref{3-5}, we obtain
\begin{align*}
N^{k+2}(a_\la b)-N^{k}(N(a)_\la N(b))=N(N^{k+1}(a)_\la b)-N^{k+1}(a)_\la N(b)+N(a_\la N^{k+1}(b))-N(a)_\la N^{k+1}(b),
\end{align*}
which can be rewritten in the following form
\begin{align*}
N^{k}(N(a))_\la N(b)+N(a)_\la N^{k}(N(b))-N^{k}(N(a)_\la N(b))
=N\big(N^{k+1}(a)_\la b+a_\la N^{k+1}(b)-N^{k+1}(a_\la b)\big).
\end{align*}
This is exactly \eqref{3-2}. Finally, applying \eqref{3-2} inductively
\begin{align*}
N^r(a\circ^{N^{k+r}}_\la b)=N^{r-1}N(a\circ^{N^{k+r}}_\la b)=N^{r-1}\big(N(a)\circ^{N^{k+r-1}}_\la N(b)\big)=N^{r-1}(N(a))\circ^{N^{k}}_\la N^{r-1}(N(b))=N^r(a)\circ^{N^k}_\la N^r(b).
\end{align*}
This ends the proof.
\end{proof}

\begin{theo}\label{th5}If $N$ is a Nijenhuis operator of $(\A,\theta_\la)$, then for arbitrary $a,b\in\A$ and $i,k=0,1,2,\cdots$, there holds
\begin{align}\label{4-1}
N^k(a)\circ^{N^i}_\la b+a\circ^{N^i}_\la N^k(b)-N^k(a\circ^{N^i}_\la b)=a\circ^{N^{i+k}}_\la b,
\end{align}
and $N^k$ is a Nijenhuis operator on $(\A, \theta^{N^i}_\la)$. In particular, all $\la$-products $\theta^{N^k}_\la$ are associative and compatible.
\end{theo}
\begin{proof} %First, we show that
%\begin{align}\label{4-2}
%(m^{N^i})^{N}_\la=m^{N^{i+1}}_\la.
%\end{align}
By a straightforward computation, we have
\begin{align*}
N^k(a)\circ^{N^i}_\la b&+a\circ^{N^i}_\la N^k(b)-N^k(a\circ^{N^i}_\la b)-a\circ^{N^{i+k}}_\la b\\
=&\cancel{N^{i+k}(a)_\la b}+N^k(a)_\la N^i(b)-N^i(N^k(a)_\la b)+N^i(a)_\la N^k(b)+\cancel{a_\la N^{i+k}(b)}-N^i(a_\la N^k(b))\\
&-N^k\big( N^i(a)_\la b+a_\la N^i(b)-N^i(a_\la b)\big)-\cancel{N^{i+k}(a)_\la b}-\cancel{a_\la N^{i+k}(b)}+N^{i+k}(a_\la b)\\
=&\big(N^i(a)_\la N^k(b)-N^k(N^i(a)_\la b)-N^i(a_\la N^k(b))+N^{i+k}(a_\la b)\big)\\
&+\big(N^k(a)_\la N^i(b)-N^i(N^k(a)_\la b)-N^k(a_\la N^i(b))+N^{i+k}(a_\la b)\big),
\end{align*}
which vanishes due to Lemma \ref{lemm3}. This proves \eqref{4-1}.

Now we apply $N^k$ to the both sides of \eqref{4-1} and obtain
\begin{align*}
N^k\big(N^k(a)\circ^{N^i}_\la b+a\circ^{N^i}_\la N^k(b)-N^k(a\circ^{N^i}_\la b)\big)=N^k(a\circ^{N^{i+k}}_\la b)=N^k(a)\circ^{N^{i}}_\la N^k(b),
\end{align*}
where we have used Lemma \ref{lemm2}. This ends the proof.
%\begin{align*}
%a(\ast^{N^i})^N_\la b&=N(a)\ast^{N^i}_\la b+a\ast^{N^i}_\la N(b)-N(a\ast^{N^i}_\la b)\\
%&=N^{i+1}(a)_\la b+N(a)_\la N^i(b)-N^i(N(a)_\la b)+N^i(a)_\la N(b)+a_\la N^{i+1}(b)-N^i(a_\la N(b))-N(a)\ast^{N^{i-1}}_\la N(b)\\
%&=N^{i+1}(a)_\la b-N^i(N(a)_\la b)+a_\la N^{i+1}(b)-N^i(a_\la N(b))-N^{i-1}(N(a)_\la N(b))\\
%&=N^{i+1}(a)_\la b+a_\la N^{i+1}(b)-N^{i+1}(a_\la b)-N^{i-1}(N(N(a)_\la b+a_\la N(b)-N(a_\la b))-N(a)_\la N(b))\\
%&=N^{i+1}(a)_\la b+a_\la N^{i+1}(b)-N^{i+1}(a_\la b)=a\ast^{N^{i+1}}_\la b,
%\end{align*}
%where we have used Lemma \ref{lemm2}. This proves the theorem.
%$N(a\ast^{N^i}_\la b)=N(a)\ast^{N^{i-1}}_\la N(b)$. Now, \eqref{4-2} and \eqref{3-2} show that
%\begin{align*}
%N(a(\ast^{N^i})^N_\la b)=N(a\ast^{N^{i+1}}_\la b)=N(a)\ast^{N^{i}}_\la N(b),
%\end{align*}
%and thus $N$ is a Nijenhuis operator on $(\A, m^{N^i}_\la)$. Finally, we can apply Lemma \ref{lemm2} and \eqref{4-2} to
%$m^{N^i}_\la$ instead of $m_\la$ that proves the theorem.
\end{proof}

There is a way to obtain a new Nijenhuis operator from Two Nijenhuis
operators. Let $N_1$ and $N_2$ be two Nijenhuis
operators on an associative conformal algebra $\A$. They are said to be {\it compatible} if $N_1+N_2$ is again a
Nijenhuis operator. Evidently, this requirement is equivalent to the following condition
\begin{eqnarray}\label{4-3}
N_1(a)_\la N_2(b)+N_2(a)_\la N_1(b)
=N_1(a\circ^{N_2}_\la b)+N_2(a\circ^{N_1}_\la b)
\end{eqnarray}
for all $a,b\in\A.$

Note that \eqref{4-3} depends linearly on $N_1$ and $N_2$. Hence, if $N_1,\cdots, N_k$ are
Nijenhuis operators on $\A$ and $N_1$ is compatible with $N_2,\cdots, N_k$, then it is compatible with any linear combination of them. If $N_1,\cdots, N_k$ are pairwise compatible, then any two linear combinations of them are compatible. %Now let $N$ be a Nijenhuis operator on $\A$. It follows from  Theorem \ref{th5} that $N^k$ for $k=0,1,2,\cdots$ are
%Nijenhuis operators on $\A$, and we have the following theorem.
\begin{theo}\label{th6}
If $N$ is a Nijenhuis operator on $\A$, then all linear combinations of $N^k$, $k=0,1,2,\cdots$, are compatible.
\end{theo}
\begin{proof} It follows from Proposition \ref{p14} that all $N^k$ for $k=0,1,2,\cdots$ are
Nijenhuis operators on $\A$.
For $k\geq r$, we have
\begin{align*}
N^k(a\circ ^{N^r}_\la b)+N^r(a\circ^{N^k}_\la b)&=N^{k-r}N^r(a\circ ^{N^r}_\la b)+N^r(a\circ_\la^{N^{k-r+r}} b)\\&=N^{k-r}(N^r(a)_\la N^r(b))+N^r(a)\circ^{N^{k-r}}_\la N^r(b)\\&=N^k(a)_\la N^r(b)+N^r(a)_\la N^k(b),
\end{align*}
where $a,b\in\A$ and we have used Lemma \ref{lemm2}. Then we get the result by \eqref{4-3}.
\end{proof}

Recall that we have introduced compatible
$\O$-operators at the end of Section 3. Here we show that
there is a close interrelation between a Nijenhuis operator and a pair of compatible
$\O$-operators.

\begin{theo}\label{th4} Let $T_1, T_2:M\rightarrow \A$ be two $\O$-operators on an associative conformal algebra $\A$ with respect to a conformal $\A$-bimodule $M$. If $T_1$ and $T_2$
are invertible, then $T_1$ and $T_2$ are compatible if and only if $N = T_1T_2^{-1}$
is a Nijenhuis operator on the associative conformal algebra $\A$.
\end{theo}
\begin{proof} For arbitrary $a,b\in\A$, there exist unique $m,n\in M$ such that $a=T_2(m)$ and $b=T_2(n)$, since $T_2$
is invertible. Hence that $N=T_1T_2^{-1}$ is a Nijenhuis operator of $\A$ is equivalent to
\begin{eqnarray}\label{5-7}
NT_2(m)_\la NT_2(n)=N\big(NT_2(m)_\la T_2(n)+T_2(m)_\la NT_2(n)-N(T_2(m)_\la T_2(n))\big).
\end{eqnarray}
As $T_1$ and $T_2$ are $\O$-operators, and $T_1=NT_2$,  \eqref{5-7} is equivalent to
\begin{eqnarray}\label{5-8}
NT_2(T_1(m)_\la n+m_\la T_1(n))=N\big(T_1(m)_\la T_2(n)+T_2(m)_\la T_1(n)-T_1(T_2(m)_\la n+m_\la T_2(n))\big).
\end{eqnarray}
Since $N$ is invertible, \eqref{5-8} is equivalent to
\begin{eqnarray*}
T_2(T_1(m)_\la n+m_\la T_1(n))=T_1(m)_\la T_2(n)+T_2(m)_\la T_1(n)-T_1(T_2(m)_\la n+m_\la T_2(n)),
\end{eqnarray*}
which is exactly \eqref{3-7}. This ends the proof.
\end{proof}

The following is a straightforward corollary of Theorems \ref{th3} and \ref{th4}.
\begin{coro} Let $T:M\rightarrow \A$ be an $\O$-operator on an associative conformal algebra $\A$ with respect to a conformal $\A$-bimodule $M$. If there exists an invertible Nijenhuis
operator $N$ on $\A$ such that $NT:M\rightarrow \A$ is also an $\O$-operator, then $(M, \succ^{T}_\la , \prec^{T}_\la)$ and $(M, \succ^{NT}_\la , \prec^{NT}_\la)$ are compatible dendriform conformal algebras, where
\begin{align*}
m\succ^{T}_\la  n=T(m)_\lambda n, ~~ m\prec^{T}_\la n=m_\lambda T(n),~~  m\succ^{NT}_\la  n=NT(m)_\lambda n, ~~ m\prec^{NT}_\la n=m_\lambda NT(n),
\end{align*}
for all $m,n\in M.$ If, in addition, $T$ is invertible, then $(\A, \succ^1_\lambda,\prec^1_\lambda)$ and
$(\A, \succ^2_\lambda,\prec^2_\lambda)$ are compatible dendriform conformal algebras, where
\begin{align*}
a\succ^1_\lambda b= T(a_\la T^{-1}(b)),~~  a \prec^1_\lambda b = T(T^{-1}(a)_\la b), ~~  a\succ^2_\lambda b= NT(a_\la (NT)^{-1}(b)),~~  a \prec^2_\lambda b = NT((NT)^{-1}(a)_\la b),
\end{align*}
for all $a,b\in \A$.

\end{coro}
\begin{prop} If $\A=\A_1\bigoplus\A_2$ is a matching pair of
associative conformal algebras $\A_1$ and $\A_2$, $P_1$ and $P_2$ denotes the corresponding projections of $\A$ onto $\A_1$ and $\A_2$, respectively, then any linear combination of $P_1$ and $P_2$ is a Nijenhuis operator of $\A$.
\end{prop}
\begin{proof} Assume that $k_1P_1+k_2P_2$ is an arbitrary linear combination of $P_1$ and $P_2$. Since $k_1P_1+k_2P_2=(k_1-k_2)P_1+k_2{\rm id}$, it is sufficient to show that $P_1$ is a Nijenhuis operator. For any $a=(a_1,a_2), b=(b_1,b_2)\in\A,$ where $a_1,b_1\in \A_1$, $a_2,b_2\in\A_2$, we have
\begin{align*}
a\circ^{P_1}_\la b&=P_1(a)_\la b+a_\la P_1(b)-P_1(a_\la b)\\
&=(a_1,0)_\la (b_1,b_2)+(a_1,a_2)_\la (b_1,0)-P_1((a_1,a_2)_\la (b_1,b_2))\\
&=(a_{1\,\la}b_1+a_1\cdot_\la^2 b_2,a_1\cdot^1_\la b_2)+(a_{1\,\la}b_1+a_2\cdot_\la^2 b_1,a_2\cdot^1_\la b_1)-(a_{1\,\la}b_1+a_2\cdot_\la^2 b_1+a_1\cdot^2_\la b_2,0)\\
&=(a_{1\,\la}b_1,a_1\cdot^1_\la b_2+a_2\cdot^1_\la b_1).
\end{align*}
Hence $P_1(a\circ^{P_1}_\la b)=(a_{1_\la}b_1,0)=P_1(a)_\la P_1(b)$.  This proves that  $P_1$ is a Nijenhuis operator.
\end{proof}

The Nihenhuis operators for Lie algebras have been widely studied. In the following, we extend this notion to Lie conformal algebras.

\begin{defi} {\rm Let $(\L,[\cdot_\la\cdot])$ be a Lie conformal algebra. A $\C[\partial]$-module homomorphism $N:\L\rightarrow \L$ is called a Nijenhuis operator of $\L$ if the condition
\begin{eqnarray}\label{NN}
[N(a)_\lambda N(b)]=N\big([N(a)_\lambda b]+[a_\lambda N(b)]-N([a_\lambda b])\big)
\end{eqnarray} is satisfied for all $a,b\in \L.$
}\end{defi}

Let $N:\L\rightarrow \L$ be a $\C[\partial]$-module homomorphism on a Lie conformal algebra $\L$. Define
\begin{align}\label{N4}
[a_\la b]^N:=[N(a)_\lambda b]+[a_\lambda N(b)]-N([a_\lambda b]), \ \forall \ a,b\in \L.
\end{align}
If $N$ is a Nijenhuis operator of $\L$, then it is not difficult to check that $(\L, [\cdot_\la\cdot]^N)$ forms a new Lie conformal algebra, and $N$ is an algebra homomorphism from $(\L, [\cdot_\la\cdot]^N)$ to $(\L,[\cdot_\la\cdot])$.

\begin{theo} If $N:\A\rightarrow \A$ is a Nijenhuis operator on an associative conformal algebra $\A$, then $N$ is also a Nijenhuis operator of the commutator Lie conformal algebra $\A^L$, and it holds
\begin{eqnarray}\label{N5}
[a_\la b]^N=a\circ_\la^N b-b\circ_{-\la-\partial}^N a,
\end{eqnarray}
for all $a,b\in\A$, namely, the deformed $\la$-bracket $[\cdot_\la\cdot]^N $ is the commutator of the deformed associative $\la$-product $\circ_\la^N$.
\end{theo}
\begin{proof}
For all $a,b\in\A$, we have
\begin{align*}
[a_\la b]^N&\stackrel{\eqref{N4}}{=}[N(a)_\lambda b]+[a_\lambda N(b)]-N([a_\lambda b])\\
&\stackrel{\eqref{L11}}{=}N(a)_\lambda b-b_{-\la-\partial}N(a)+a_\lambda N(b)-N(b)_{-\la-\partial}a-N(a_\lambda b-b_{-\la-\partial}a)\\
&=\big(N(a)_\lambda b+a_\lambda N(b)-N(a_\lambda b)\big)-\big( N(b)_{-\la-\partial}a+ b_{-\la-\partial}N(a)-N(b_{-\la-\partial}a) \big)\\&\stackrel{\eqref{ast}}{=}a\circ_\la^N b-b\circ_{-\la-\partial}^N a.
\end{align*}
Then it follows that
\begin{align*}
N([a_\la b]^N)\stackrel{\eqref{N5}}{=}N\big(a\circ_\la^N b-b\circ_{-\la-\partial}^N a\big)\stackrel{\eqref{N3}}{=}N(a)_\la N(b)-N(b)_{-\la-\partial}N(a)\stackrel{\eqref{L11}}{=}[N(a)_\la N(b)].
\end{align*}
This shows that $N$ is a Nijenhuis operator of $\A^L$.
\end{proof}

More properties of Nijenhuis operators of Lie conformal algebras will be given in our next paper \cite{YZ}.

\subsection{Deformations of associative conformal algebras}

Let $\A$ be an associative conformal algebra, and $\omega_\la:\A\times\A\rightarrow\A[\la]$ is a conformal bilinear map. We consider a $t$-parameterized family of bilinear $\la$-multiplications
\begin{align}\label{4-4}
a\circ^t_\la b=a_\la b+t \omega_\la(a,b),
\end{align}
where $a,b\in\A.$ If all the $\la$-multiplications $\circ^t_\la$ endow $\A$ with associative conformal algebra structures, then we say that $\omega$ generates a {\it deformation} of the associative conformal algebra $\A$. Evidently, this requirement is
equivalent to the conditions
\begin{align}
\omega_{\la+\mu}(a_\la b,c)+\omega_\la(a,b)_{\la+\mu}c&=\omega_\la(a,b_\mu c)+a_\la\omega_\mu(b,c),\label{4-5}\\
\omega_{\la+\mu}(\omega_\la(a,b),c)&=\omega_\la(a,\omega_\mu(b,c)).\label{4-6}
\end{align}
Hence, $\omega_\la$ must itself be an associative conformal algebra structure (cf. \eqref{4-6}), satisfying condition \eqref{4-5}. Recalling the definition of the coboundary operator in the Hochschild cohomology complex of $\A$
with the adjoint action of $\A$, we can present \eqref{4-5} in an
abbreviated form, $\mathbf{d}(\omega)=0$. Namely, $\omega$ is a $2$-cocyle in $C^2(\A,\A)$.

A deformation $\omega$ is said to be {\it trivial} if there exists a $\C[\partial]$-module homomorphism $N: \A\rightarrow \A$ such that for $T_t={\rm id}+t N$ there holds
\begin{align}\label{4-7}
T_t(a\circ^t_\la b)=T_t(a)_\la T_t(b), \ \text{for\ all} \ a,b\in\A.
\end{align}
As we have
\begin{align*}
T_t(a\circ^t_\la b)=a_\la b+t(\omega_\la(a,b)+N(a_\la b))+t^2 N\omega_\la(a,b),
\end{align*}
and \begin{align*}
 T_t(a)_\la T_t(b)=a_\la b+t(N(a)_\la b+a_\la N(b))+t^2 N(a)_\la N(b),
\end{align*}
the triviality of deformation is equivalent to the conditions
\begin{align}
\omega_\la(a,b)&=N(a)_\la b+a_\la N(b)-N(a_\la b),\label{4-8}\\
 N\omega_\la(a,b)&= N(a)_\la N(b).\label{4-9}
\end{align}
It follows from \eqref{4-8} and \eqref{4-9} that $N$ must satisfy the following condition:
\begin{align*}
 N(a)_\la N(b)=N\big(N(a)_\la b+a_\la N(b)-N(a_\la b)\big),
\end{align*}
which is to say that $N$ is a Nijenhuis operator of $\A$.

We have deduced that any trivial deformation produces a Nijenhuis
operator. Notably, the converse is also valid, as the following theorem shows.
\begin{theo} Let $N:\A\rightarrow \A$ ba a Nijenhuis operator. Then a deformation of $\A$
can be obtained by putting
\begin{align}
 \omega_\la(a,b)=N(a)_\la b+a_\la N(b)-N(a_\la b), \label{4-10}
\end{align}
for all $a,b\in\A$. This deformation is a trivial one.
\end{theo}
\begin{proof} By Corollary \ref{coro1}, we have $\omega_\la$ is associative, namely, \eqref{4-6} is valid. As \eqref{4-10} can be represented in terms of the coboundary
operator  in the Hochschild cohomology complex of $\A$
with the adjoint action of $\A$ as $\omega=\mathbf{d}(N)$, we have $\mathbf{d}(\omega)=0$ and therefore
condition \eqref{4-5} holds.

Evidently, \eqref{4-8} and \eqref{4-9} are satisfied and therefore $\omega$ generates a trivial deformation of $\A$.  This ends the proof.
\end{proof}

\section{Cohomology of $\O$-operators}

In this section, we study the cohomology problem of $\O$-operators by using Gerstenhaber-bracket \cite{Ger} and derived bracket construction of Kosmann-Schwarzbach \cite{KS1,KS2} (see also \cite{Das,U1,U2,U3}).

%\subsection{Cohomology}

Let $T:M\rightarrow \A$ be an $\mathcal{O}$-operator on an associative conformal algebra $\A$ with respect to a conformal bimodule $M$. Then $M$ carries an associative conformal algebra structure  $ M_{ass}$ given by \eqref{M-ass} and there is a conformal $M_{ass}$-bimodule structure on $\A$ given by Lemma \ref{lemm1}. The Hochschild cohomology of $M_{ass}$ with coefficients in the conformal bimodule $\A$ is by definition the cohomology of the $\mathcal{O}$-operator $T$. More precisely, the space of $n$-cochains
$C^n(M,\A)$ for $n\geq 1$ consists of all multilinear maps of the form
$$f_{\lambda_1,\cdots,\lambda_{n-1}}:M^{\otimes n}\longrightarrow \A[\lambda_1, \cdots , \lambda_{n-1}]$$
satisfying the sesquilinearity conditions \eqref{co1} and \eqref{co2}. It follows from \eqref{co3}, \eqref{M-ass}, and Lemma \ref{lemm1} that
%\begin{align}
%f_{\lambda_1,\cdots,\lambda_{n-1}}(u_1,\cdots,\partial u_i,\cdots,u_n)&=-\la_i f_{\lambda_1,\cdots,\lambda_{n-1}}(u_1,\cdots,u_n), ~~i=1,\cdots,n-1,\\
%f_{\lambda_1,\cdots,\lambda_{n-1}}(u_1,\cdots,u_{n-1},\partial u_n)&=(\partial+\lambda_1+\cdots+\lambda_{n-1})f_{\lambda_1,\cdots,\lambda_{n-1}}(u_1,\cdots,u_n).
%\end{align}
the conformal Hochschild differential $\mathbf{d}:C^n(M,\A)\rightarrow C^{n+1}(M,\A)$ is given by
\begin{align}\label{6-1-1}
(\mathbf{d}(f))_{\lambda_1,\cdots,\lambda_{n}}(u_1,\cdots,u_{n+1})=&T(u_1)_{\la_1}f_{\lambda_2,\cdots,\lambda_n}(u_2,\cdots,u_{n+1})-T\big(u_{1\,{\la_1}}f_{\lambda_2,\cdots,\lambda_n}(u_2,\cdots,u_{n+1})\big)\nonumber\\
&+\sum_{i=1}^{n}(-1)^i f_{\lambda_1,\cdots,\lambda_i+\lambda_{i+1},\cdots,\lambda_{n}}(u_1,\cdots,T(u_i)_{\la_i}u_{i+1}+u_{i\,{\la_i}\,}T(u_{i+1}),\cdots,u_{n+1})\nonumber\\
&+(-1)^{n+1}f_{\lambda_1,\cdots,\lambda_{n-1}}(u_1,\cdots,u_n)_{\la_1+\cdots+\la_n}T(u_{n+1})\nonumber\\
&-(-1)^{n+1}T\big(f_{\lambda_1,\cdots,\lambda_{n-1}}(u_1,\cdots,u_n)_{\la_1+\cdots+\la_n}u_{n+1}\big).
\end{align}

%We set $C^0(M,\A)=\A$. For $a\in\A$, we define
%\begin{align}
%\mathbf{d}(a)(m):=m\cdot_\lambda a-a\cdot_\lambda m=T(m)_\lambda a-T(m_\lambda a)-a_\lambda T(m)+T(a_\lambda m),\label{5-6}
%\end{align} for all $m\in M.$

Denote by $Z^n(M,\A)$ and $B^n(M,\A)$ the spaces of $n$-cocycles and $n$-coboundaries,
respectively. Then the quotient space $$H^n(M,\A) = Z^n(M,\A)/B^n
( M,\A)$$ is called the $n$th Hochschild cohomology group of $M_{ass}$ with coefficients in $\A$.

For instance, the space of $1$-cocycles $Z^1( M,\A) = \Ker \,{\mathbf{d}} \subseteq  C^1(M,\A)$ consists of all
$\C[\partial]$-linear maps $f: M\rightarrow \A$ such that
\begin{align}\label{1-cocycle-1}
0=(\mathbf{d}(f))_{\lambda}(u,v)=T(u)_\la f(v)+f(u)_\la T(v)-T\big(u_\la f(v)+f(u)_\la v\big)-f\big(u_\la T(v)+T(u)_\la v\big),
\end{align}
for all $u,v\in M.$ Then the following result is straightforward.
\begin{prop}\label{p15}
$T\in C^1(M,\A)$ is an $\O$-operator if and only if ${\mathbf{d}}(T)=0$.
\end{prop}

 In the literature \cite{Wu}, Wu generalized the Gerstenhaber-bracket \cite{Ger} to the pseudotensor category and constructed a differential graded
Lie algebra which controls the cohomology theory of $H$-pseudoalgebras. As associative conformal algebras
form a special class of $H$-pseudoalgebras, we can translate the construction of Wu in the case of associative conformal algebras by means of $\la$-products as follows.

Let $(\A,\theta_\la)$ be an associative conformal algebra with $\theta_\la(a,b):=a_\la b$ for all $a,b\in\A$. For $n\geq 1$, set $C^n(\A)=C^{n}(\A,\A)$ and $C^\bullet(\A)=\bigoplus_{n\geq 1} C^n(\A)$. For any $f\in C^{m}(\A)$ and $g\in C^{n}(\A)$, define the G-bracket on $C^\bullet(\A)$ by
\begin{align}\label{G}
[f,g]:=f\circ g-(-1)^{(m-1)(n-1)}g\circ f,
\end{align}
where $f\circ g\in C^{m+n-1}(\A)$ is defined by
\begin{align*}
(f\circ g)&_{\la_1,\cdots,\la_{m+n-2}}(a_1,\cdots,a_{m+n-1})\\
=&\sum_{i=1}^{m}(-1)^{(i-1)(n-1)}f_{\la_1,\cdots,\la_{i-1},\la_i+\cdots+\la_{i+n-1},\la_{i+n},\cdots,\la_{m+n-2}}(a_1,\cdots,a_{i-1},g_{\la_i,\cdots,\la_{i+n-2}}(a_i,\cdots,a_{i+n-1}),\cdots,a_{m+n-1}),
\end{align*}
for all $a_1,\cdots,a_{m+n-1}\in\A.$

\begin{lemm}\label{lemm7} (\cite[Lemma 3.2]{Wu}) Let $S\in C^{2}(\A)$. Then $[S,S]=0$ if and
only if $S$ is associative.
\end{lemm}
 It follows that $[\theta,\theta]=0$. For any $f\in C^{n}(\A)$, we have
\begin{align*}
[\theta,f]&_{\la_1,\cdots,\la_n}(a_1,\cdots,a_{n+1})\\&=\big(\theta \circ f-(-1)^{n-1}f\circ \theta\big)_{\la_1,\cdots,\la_{n}}(a_1,\cdots,a_{n+1})\\
&=\theta_{\la_1+\cdots+\la_{n}}(f_{\la_1,\cdots,\la_{n-1}}(a_1,\cdots,a_{n}),a_{n+1})+(-1)^{n-1}\theta_{\la_1}(a_1,f_{\la_2,\cdots,\la_{n}}(a_2,\cdots,a_{n+1}))\\
&~~~~~ -(-1)^{n-1}\sum_{i=1}^{n}(-1)^{i-1}f_{\la_1,\cdots,\la_{i-1},\la_i+\la_{i+1},\la_{i+2},\cdots,\la_{n}}(a_1,\cdots,a_{i-1},\theta_{\la_i}(a_i,a_{i+1}),a_{i+2},\cdots,a_{n+1})
\\&=f_{\la_1,\cdots,\la_{n-1}}(a_1,\cdots,a_{n})_{\la_1+\cdots+\la_{n}}a_{n+1}+(-1)^{n-1}a_{1\,\la_1}f_{\la_2,\cdots,\la_{n}}(a_2,\cdots,a_{n+1})\\
&~~~~~+(-1)^{n-1}\sum_{i=1}^{n}(-1)^{i}f_{\la_1,\cdots,\la_{i-1},\la_i+\la_{i+1},\la_{i+2},\cdots,\la_{n}}(a_1,\cdots,a_{i-1},a_{i\,\la_i}a_{i+1},a_{i+2},\cdots,a_{n+1}).
\end{align*}
Comparing this with \eqref{co3}, we obtain $[\theta,f]=(-1)^{n-1}\mathbf{d}(f)$.  By
the graded Jacobi identity of G-bracket, $d_\theta:=[\theta,\cdot]$ becomes a square-zero derivation
of degree $+1$.

\begin{theo}\label{th8} (\cite[Theorem 3.3]{Wu}) If $(\A,\theta_\la)$ is an associative conformal algebra, then $C^\bullet(\A)$ is a dg-Lie algebra with the G-bracket defined by \eqref{G} and the differential $d_\theta:=[\theta,\cdot]$.
\end{theo}

Let $M$ be a conformal bimodule over the associative conformal algebra $(\A,\theta_\la)$. In the following, we consider the G-bracket defined by \eqref{G} on the graded vector space $C^\bullet(\A\oplus _0 M)=\bigoplus_{n\geq 1}C^{n}(\A\oplus _0 M,\A\oplus _0 M)$. Denote by $\hat{\theta}_\la$ the associative $\la$-multiplication of $\A\oplus _0 M$, i.e.,
\begin{align}\label{hat-theta}
\hat{\theta}_\la (a,b)=a_\la b,
\ \hat{\theta}_\la (a,n)=a_\la n,\ \hat{\theta}_\la (m,b)=m_\la b, \ \hat{\theta}_\la (m,n)=0,
\end{align}
for all $a,b\in\A$ and $m,n\in M$. By Lemma \ref{lemm7}, we have $[\hat{\theta},\hat{\theta}]=0$. By
the graded Jacobi identity of G-bracket, $d_{\hat\theta}:=[{\hat\theta},\cdot]$ becomes a square-zero derivation
of degree $+1$. By using this derivation, we define a derived bracket (cf. \cite{KS1,KS2}) on $C^\bullet(\A\oplus _0 M)$ by
\begin{align}\label{derived bracket}
[[f,g]]:=(-1)^m[d_{\hat\theta}(f),g]=(-1)^m[[\hat\theta,f],g],
\end{align}
where $f\in C^m(\A\oplus _0 M)$ and $g\in C^n(\A\oplus _0 M)$. It is worth noticing that the derived bracket $[[\cdot,\cdot]]$ is not
graded commutative, but it satisfies the graded Leibniz rule (cf. \eqref{1-3}).

Recall that we have defined a lift $\hat{T}:\A\oplus M\rightarrow \A\oplus M$ for any $\C[\partial]$-module homomorphism $T:M\rightarrow \A$ by $\hat{T}(a,m):=(T(m),0)$, for all $a\in\A$ and $m \in M.$
Then the derivation of $\hat{T}$ by $\hat{\theta}$ has the the following form
\begin{align}\label{6-7}
[\hat{\theta},\hat{T}]_\la=\hat{\theta}_\la(\hat{T}\otimes{\rm id}+{\rm id}\otimes\hat{T})-(\hat{T}\circ \hat{\theta})_\la.
\end{align}
%where $\hat{T}$ is defined in Proposition \ref{p3-3}.
\begin{prop}\label{pp17} $T:M\rightarrow \A$ is an  $\mathcal{O}$-operator if and only if $\hat{T}$ satisfies $[[\hat{T},\hat{T}]]=0$.
%is a solution of Maurer-Cartan equation, i.e., $[[\hat{T},\hat{T}]]=0$.
\end{prop}
\begin{proof}For any $a,b\in\A$ and $m,n\in M$, we have $\hat{T}\circ \hat{T}(a,m)=\hat{T}(\hat{T}(a,m))=\hat{T}(T(m),0)=(T(0),0)=0$. With this and \eqref{6-7}, we have
\begin{align*}
[[\hat{T},\hat{T}]]_\la((a,m),(b,n))%=-[[\hat\theta,\hat{T}],\hat{T}]_\la((a,m),(b,n))
=&-([\hat\theta,\hat{T}]\circ\hat{T}-\hat{T}\circ[\hat\theta,\hat{T}])_\la((a,m),(b,n))\\
%=&[\hat\theta,\hat{T}]_\la(\hat{T}(a,m),(b,n))+[\hat\theta,\hat{T}]_\la((a,m),\hat{T}(b,n))-\hat{T}\circ \hat\theta_\la(\hat{T}(a,m),(b,n))-\hat{T}\circ \hat\theta_\la((a,m),\hat{T}(b,n))\\
=&-[\hat\theta,\hat{T}]_\la((T(m),0),(b,n))-[\hat\theta,\hat{T}]_\la((a,m),(T(n),0))\\&+\hat{T} \hat\theta_\la((T(m),0),(b,n))+\hat{T} \hat\theta_\la((a,m),(T(n),0))\\
%=&\hat{\theta}_\la(\hat{T}\otimes{\rm id}+{\rm id}\otimes\hat{T})_\la((T(m),0),(b,n))+\hat{\theta}_\la(\hat{T}\otimes{\rm id}+{\rm id}\otimes\hat{T})_\la((a,m),(T(n),0))
%\\&-2\hat{T}\circ \hat\theta_\la((T(m),0),(b,n))-2\hat{T}\circ \hat\theta_\la((a,m),(T(n),0))\\
=&-2\hat{\theta}_\la((T(m),0),(T(n),0))+2\hat{T}(T(m)_\la b,T(m)_\la n)+2\hat{T}(a_\la T(n),m_\la T(n))\\
=&-2\big(T(m)_\la T(n)-T(T(m)_\la n+m_\la T(n)),0\big).
\end{align*}
This proves the result.
\end{proof}

From the graded Jacobi rule of G-bracket, we obtain the following corollary.
\begin{coro}\label{coro3} If $T:M\rightarrow \A$ is an  $\mathcal{O}$-operator, then  $[\hat{\theta},\hat{T}]_\la$ is  associative on $\A\oplus _0 M$ of the form
\begin{align*}
[\hat\theta,\hat{T}]_\la((a,m),(b,n))=(m\cdot_\la b+a\cdot_\la n, m\star_\la n),
\end{align*}
for all $a,b\in\A$ and $m,n\in M$, where $\star_\la$ is defined by \eqref{M-ass} and $\cdot_\la$ is defined by \eqref{2-1}.
\end{coro}
\begin{proof} Since $[[\hat\theta,\hat{T}],[\hat\theta,\hat{T}]]=0$, %$[[[\hat\theta,\hat{T}],\hat\theta],\hat{T}]-\underbrace{[\hat\theta,[[\hat\theta,\hat{T}],\hat{T}]]}_{=0}=[[\hat\theta,[\hat\theta,\hat{T}]],\hat{T}]=0,$
$[\hat{\theta},\hat{T}]$ is associative. By \eqref{6-7}, we have
\begin{align*}
[\hat\theta,\hat{T}]_\la((a,m),(b,n))&=(T(m)_\la b-T(m_\la b)+a_\la T(n)-T(a_\la n),T(m)_\la n+m_\la T(n))\\
&=(m\cdot_\la b+a\cdot_\la n, m\star_\la n),
\end{align*}
for all $a,b\in\A$ and $m,n\in M.$  This ends the proof.
\end{proof}
\begin{remark}{\rm The associativity of $[\hat{\theta},\hat{T}]$ in  Corollary \ref{coro3} also implies the bimodule action of $M_{ass}$ on $\A$. This gives a second proof of Lemma \ref{lemm1}.}
\end{remark}

In the following, we consider the graded subspace $C^\bullet(M, \A)=\bigoplus_{n\geq 1} C^{n}(M, \A)$ of $C^\bullet(\A\oplus_0 M)$. Let $f$ be a $k$-cochain in $C^{k}(M, \A)$. We can construct a $k$-cochain $\hat{f}\in C^{k}(\A\oplus _0 M)$ by
\begin{align}
\hat{f}_{\la_1,\cdots,\la_{k-1}}((a_1,u_1), \cdots, (a_k,u_k))=(f_{\la_1,\cdots,\la_{k-1}}(u_1,\cdots,u_k),0)
\end{align}
for all $a_1,\cdots,a_k\in\A$ and $u_1,\cdots, u_k\in M.$ We call the cochain $\hat f$ a {\it horizontal lift} or simply {\it lift} of $f$. Then the graded
space $C^\bullet(M, \A)$ is identified with an abelian subalgebra of the differential graded Lie algebra
$C^\bullet(\A\oplus _0 M)$
via the horizontal lift.
One can easily check that the derived bracket defined by \eqref{derived bracket}
is closed in $C^\bullet(M, \A)$ , %and $[[T,\cdot]]$ is a square-zero derivation of degree $+1$ for the derived bracket.
Hence, by Lemma \ref{lemm55}, $(C^\bullet(M, \A), [[\cdot,\cdot]])$ becomes a graded Lie algebra. More precisely, the Lie bracket on $(C^\bullet(M, \A)$ is given by
\begin{align}\label{666}
&[[f,g]]_{\la_1,\cdots,\la_{m+n-1}}(u_1,\cdots,u_{m+n})\nonumber\\
=&(-1)^{mn}\big\{f_{\la_1,\cdots,\la_{m-1}}(u_1,\cdots,u_m)_{\la_1+\cdots+\la_m} g_{\la_{m+1},\cdots,\la_{m+n-1}}(u_{m+1},\cdots,u_{m+n})\nonumber\\ &-(-1)^{mn}g_{\la_1,\cdots,\la_{n-1}}(u_1,\cdots,u_n)_{\la_1+\cdots+\la_n}f_{\la_{n+1},\cdots,\la_{m+n-1}}(u_{n+1},\cdots,u_{m+n})
\big\}\nonumber\\
&+\big\{\sum_{i=1}^{m}(-1)^{(i-1)n} f_{\la_1,\cdots,\la_i+\cdots+\la_{i+n},\cdots,\la_{m+n-1}}(u_1,\cdots,g_{\la_i,\cdots,\la_{i+n-2}}(u_{i},\cdots,u_{i+n-1})_{\la_i+\cdots+\la_{i+n-1}}u_{i+n},\cdots,u_{m+n})\nonumber\\
&~~~~ ~~ -\sum_{i=1}^{m}(-1)^{in}f_{\la_1,\cdots,\la_{i-1},\la_i+\cdots+\la_{i+n},\cdots,\la_{m+n-1}}(u_1,\cdots,u_{i-1},u_{i\,\la_i}g_{\la_{i+1},\cdots,\la_{i+n-1}}(u_{i+1},\cdots,u_{i+n}),\cdots,u_{m+n})\big\}\nonumber\\
&-(-1)^{mn}\big\{\sum_{i=1}^{n}(-1)^{(i-1)m} g_{\la_1,\cdots,\la_i+\cdots+\la_{i+m},\cdots,\la_{m+n-1}}(u_1,\cdots,f_{\la_i,\cdots,\la_{i+m-2}}(u_{i},\cdots,u_{i+m-1})_{\la_i+\cdots+\la_{i+m-1}}u_{i+m},\cdots,u_{m+n})\nonumber\\
&~~~~ ~~-\sum_{i=1}^{n}(-1)^{im}g_{\la_1,\cdots,\la_{i-1},\la_i+\cdots+\la_{i+m},\cdots,\la_{m+n-1}}(u_1,\cdots,u_{i-1},u_{i\,\la_i}f_{\la_{i+1},\cdots,\la_{i+m-1}}(u_{i+1},\cdots,u_{i+m}),\cdots,u_{m+n})\big\},
\end{align}
where $f\in C^m(M,\A)$, $g\in C^n(M,\A)$ and $u_1,\cdots,u_{m+n}\in M.$

It follows from \eqref{666} that for arbitrary $T\in C^1(M,\A)$, $f\in C^n(M,\A)$ and  $u_1,\cdots,u_{n+1}\in M,$ there holds
\begin{align}
[[T,f]]&_{\la_1,\cdots,\la_n}(u_1,\cdots,u_{n+1})\nonumber\\=&(-1)^{n}T(u_1)_{\la_1}f_{\la_2,\cdots,\la_n}(u_2,\cdots,u_{n+1})-f_{\la_1,\cdots,\la_{n-1}}(u_1,\cdots,u_n)_{\la_1+\cdots+\la_n}T(u_{n+1})
\nonumber\\
&+T(f_{\la_1,\cdots,\la_{n-1}}(u_1,\cdots,u_n)_{\la_1+\cdots+\la_n}u_{n+1})-(-1)^{n}T(u_{1\,\la_1}f_{\la_2,\cdots,\la_n}(u_2,\cdots,u_{n+1}))\nonumber\\
&+(-1)^{n}\sum_{i=1}^{n}(-1)^{i} f_{\lambda_1,\cdots,\lambda_i+\lambda_{i+1},\cdots,\lambda_{n}}(u_1,\cdots,u_{i-1},T(u_i)_{\la_i}u_{i+1}+u_{i\,{\la_i}\,}T(u_{i+1}),u_{i+2},\cdots,u_{n+1}). \label{6-9}
\end{align}
In particular, for $T,T^\prime\in C^1(M,\A)$ and $u,v\in M$, we have
\begin{eqnarray}
 [[T,T^\prime]]_{\la}(u,v)=T(T^\prime(u)_\la v+u_\la T^\prime(v))+T^\prime(T(u)_\la v+u_\la T(v))-T(u)_\la T^\prime (v)-T^\prime(u)_\la T (v).
\end{eqnarray}
 Combining \eqref{6-9} with \eqref{6-1-1}, we obtain that if $T$ is an $\O$-operator, then
\begin{align}\label{6-8}
d_T(f):=[[T,f]]=(-1)^{n}\mathbf{d}(f).
\end{align}
Hence $d_T:=[[T,\cdot]]$ is a square-zero derivation of degree $+1$, and %$T$ is a Maurer-Cartan element in $(C^\bullet(M, \A), [[\cdot,\cdot]])$, i.e.,
$T$ satisfies $[[T,T]]=0.$

We summarize the above discussions in the following theorem.

\begin{theo}\label{th9} Let $M$ be a conformal bimodule over an associative conformal algebra $\A$.
\begin{itemize}\item[(1)]The graded vector space $C^\bullet(M, \A)=\bigoplus_{n\geq 1} C^{n}(M, \A)$ together with the bracket $[[\cdot,\cdot]]$  defined by \eqref{derived bracket} forms a graded Lie algebra. An element $T$ in $C^{1}(M, \A)$ is an $\O$-operator if and only if $T$ satisfies $[[T,T]]=0.$ %the Maurer-Cartan equation, i.e.,
\item[(2)] If $T:M\rightarrow \A$ is an $\O$-operator, then $T$ induces a differential $d_T=[[T,\cdot]]$ which makes the graded Lie algebra
$(C^\bullet(M, \A), [[\cdot,\cdot]])$ into a dg-Lie algebra. Moreover, for any $T^\prime\in C^1(M,\A)$, $T+T^\prime$ is still an $\O$-operator if and only if $T^\prime$ is a Maurer-Cartan element in $(C^\bullet(M, \A), [[\cdot,\cdot]],d_T)$, i.e., it
  satisfies $d_T(T^\prime)+\frac12[[T^\prime,T^\prime]]=0$.
\end{itemize}
\end{theo}

\begin{remark}{\rm For an $\O$-operator $T:M\rightarrow\A$, we have obtained two cochain complexes, i.e., $(C^\bullet(M,\A),\mathbf{d})$ and $(C^\bullet(M,\A), d_T)$. But the corresponding cohomologies are isomorphic by \eqref{6-8}. Hence we may use the same notation $H^\bullet(M, \A)$ to denote
the cohomology of an $\O$-operator $T$.}
\end{remark}

As we have mentioned that a Rota-Baxter operator (of weight $0$) on an associative conformal algebra
$\A$ can be seen as an $\O$-operator on $\A$ with respect to the adjoint bimodule $\A$. Therefore,
by considering the adjoint bimodule instead of arbitrary bimodule, we get a similar result of Theorem \ref{th9}.

\begin{theo} Let $\A$ be an associative conformal algebra. Then
\begin{itemize}\item[(1)] the graded vector space $C^\bullet(\A, \A)=\bigoplus_{n\geq 1} C^{n}(\A, \A)$ has a graded Lie algebra structure $[[\cdot,\cdot]]$ defined by \eqref{derived bracket}. An element $T\in C^{1}(\A, \A)$ is a Rota-Baxter operator of weight $0$ if and only if $T$ satisfies $[[T,T]]=0.$
\item[(2)] If $T\in C^{1}(\A, \A)$ is a Rota-Baxter operator of weight $0$,  then it induces a differential $d_T=[[T,\cdot]]$ on the graded Lie algebra $(C^\bullet(\A, \A), [[\cdot,\cdot]])$ to make it into a dg-Lie algebra. Further, for any $T^\prime\in C^1(\A,\A)$, $T+T^\prime$ is again a Rota-Baxter operator of weight $0$ if and only if $T^\prime$ satisfies the Maurer-Cartan equation, i.e., $d_T(T^\prime)+\frac12[[T^\prime,T^\prime]]=0$.
\end{itemize}
\end{theo}

Given a Rota-Baxter operator $T$ of weight $0$ on an associative conformal algebra $\A$, the vector space $\A$
carries a dendriform conformal algebra structure (cf. Proposition \ref{p3}). Hence, by \eqref{M-ass-1}, $\A$ carries a new associative $\la$-product $a\star_\la b = a_\la T(b)+T(a)_\la b,$  for $a, b\in \A$. Denote this associative conformal algebra by $\A_{ass}$. By Lemma \ref{lemm1}, $\A_{ass}$ has a conformal bimodule action on $\A$ given by
\begin{align}\label{2-1-1}
a\cdot_\lambda b=T(a)_\lambda b-T(a_\lambda b), ~~ b\cdot_\lambda a=b_\lambda T(a)-T(b_\lambda a),
\end{align}
for $a\in \A_{ass}$ and $b\in\A$. Then the cohomology of the associative conformal algebra $\A_{ass}$ with coefficients in the above conformal bimodule structure on $\A$ is called the {\it cohomology} of the Rota-Baxter operator $T$.

\begin{remark}{\rm  Note that the associative conformal algebra $\A_{ass}=(\A, \star_\la)$ has two more conformal bimodule structures
on $\A$. The first one is given by the adjoint bimodule $a\cdot_\lambda b=a\star_\lambda b$ and $b\cdot_\lambda a=b\star_\lambda a$. The second one is given by $a\cdot_\lambda b=T(a)_\lambda b$ and $b\cdot_\lambda a=b_\lambda T(a)$. However, neither of
these two bimodule structures are the same (in general) with that of \eqref{2-1-1}.}
\end{remark}

%Let $T:M\rightarrow \A$ be an $\mathcal{O}$-operator. Now we revisit the matching pair $\A\bowtie M_{ass}$, which has an associative $\la$-multiplication of the form \eqref{M-P}.

Finally, we consider the twisted case. Let $M$ be a conformal bimodule of an associative conformal algebra $\A$. For any $T\in C^1(M,\A)$, we can see from the proof of Proposition \ref{pp17} that
\begin{align}\label{6-10}
\frac12[[\hat{T},\hat{T}]]_\la=\hat{T}\circ\hat{\theta}_\la(\hat{T}\otimes{\rm id}+{\rm id}\otimes\hat{T})-\hat{\theta}_\la(\hat{T}\otimes \hat{T}),
\end{align}where $\hat{\theta}_\la$ is defined by \eqref{hat-theta}.
Let $\varphi_\la\in C^2(\A,M)$ be a $2$-cocycle. We define a 2-cochain $\hat{\varphi}\in C^2(\A\oplus_0 M)$ by
\begin{align*}
\hat{\varphi}_\la (a,b)={\varphi}_\la(a,b),~~ \hat{\varphi}_\la(m,b)=\hat{\varphi}_\la(a,n)=\hat{\varphi}_\la(m,n)=0,
\end{align*}
for all $a,b\in\A$ and $m,n\in M.$
\begin{lemm}\label{lemm6} $T\in C^1(M,\A)$ is a $\varphi$-Rota-Baxter operator if and only if $\hat{T}$ satisfies
\begin{align}\label{6-11}
\hat{\theta}_\la(\hat{T}\otimes \hat{T})-\hat{T}\circ\hat{\theta}_\la(\hat{T}\otimes{\rm id}+{\rm id}\otimes\hat{T})-\hat{T}\circ \hat{\varphi}_\la(\hat{T}\otimes \hat{T})=0.
\end{align}
\end{lemm}
\begin{proof} For all $a,b\in\A$ and $m,n\in M$, we have
\begin{align*}
\big(\hat{\theta}_\la(\hat{T}\otimes \hat{T})&-\hat{T}\circ\hat{\theta}_\la(\hat{T}\otimes{\rm id}+{\rm id}\otimes\hat{T})-\hat{T}\circ \hat{\varphi}_\la(\hat{T}\otimes \hat{T})\big)((a,m),(b,n))\\&=\Big(T(m)_\la T(n)-T\big(T(m)_\la n+m_\la T(n)+\varphi_\la(T(m),T(n))\big),0\Big).
\end{align*}
This implies the result.
\end{proof}

\begin{prop} $T:M\rightarrow \A$ is a $\varphi$-Rota-Baxter operator if and only if $\hat{T}$ satisfies the modified Maurer-Cartan equation:
\begin{align*}
\frac12[[\hat{T},\hat{T}]]_\la=-\frac16[[[\hat{\varphi},\hat{T}],\hat{T}],\hat{T}]_\la.
\end{align*}
\end{prop}
\begin{proof} For any $a,b\in\A$ and $m,n\in M$, we have
\begin{align*}
\frac12[[[\hat{\varphi},\hat{T}],\hat{T}],\hat{T}]_\la((a,m),(b,n)){=}&\frac12([[\hat{\varphi},\hat{T}],\hat{T}]\circ \hat{T})_\la((a,m),(b,n))-\frac12 (\hat{T}\circ [[\hat{\varphi},\hat{T}],\hat{T}])_\la((a,m),(b,n))\\
{=}&(\hat{T}\circ\hat{\varphi}_\la(\hat{T}\otimes{\rm id}+{\rm id}\otimes\hat{T})-\hat{\varphi}_\la(\hat{T}\otimes \hat{T}))((T(m),0),(b,n))\\
&+(\hat{T}\circ\hat{\varphi}_\la(\hat{T}\otimes{\rm id}+{\rm id}\otimes\hat{T})-\hat{\varphi}_\la(\hat{T}\otimes \hat{T}))((a,m),(T(n),0))\\&-\hat{T}(\hat{T}\circ\hat{\varphi}_\la(\hat{T}\otimes{\rm id}+{\rm id}\otimes\hat{T})-\hat{\varphi}_\la(\hat{T}\otimes \hat{T}))((a,m),(b,n))\\
=&3\hat{T}\circ\hat{\varphi}_\la(\hat{T}\otimes \hat{T}))((a,m),(b,n)),
\end{align*}
which gives $\hat{T}\circ\hat{\varphi}_\la(\hat{T}\otimes \hat{T}))=\frac16[[[\hat{\varphi},\hat{T}],\hat{T}],\hat{T}]_\la$. From \eqref{6-10} and Lemma \ref{lemm6}, we obtain the desired result.
\end{proof}

 \vs{10pt}
\noindent{\bf{ Acknowledgements.}}\ {Supported by National Natural Science
Foundation grants of China (11301109). }


\small
\begin{thebibliography}{9999}\vskip0pt\small{
\def\re{\bibitem}\parindent=2ex\parskip=-2pt\baselineskip=-2pt


\bibitem{A} F.V. Atkinson, Some aspects of Baxter's functional equation, J. Math. Anal. Appl. 7 (1) (1963) 1--30.

\bibitem {A2} M. Aguiar, Pre-Poisson algebras, Lett. Math. Phys. 54 (4) (2000) 263--277.

\bibitem {A3} M. Aguiar, On the associative analog of Lie bialgebras, J. Algebra 244 (2001) 492--532.

%\bibitem {A4} M. Aguiar, Infinitesimal bialgebras, pre-Lie and dendriform algebras, in: Hopf Algebras, in: Lect. Notes Pure
%Appl. Math., vol. 237, 2004, pp. 1--33.


%\bibitem {AB} H. An, C. Bai, From Rota-Baxter algebras to pre-Lie algebras, J. Phys. A 41 (2008) 015201.

\bibitem {B} G. Baxter, An analytic problem whose solution follows from a simple algebraic identity, Pacific J.
Math. 10 (3) (1960) 731--742.


%\bibitem {BGN} C. Bai, L. Guo, X. Ni, Nonabelian generalized lax pairs, the classical Yang-Baxter equation and
%Post-Lie algebras, Commun. Math. Phys. 297 (2) (2010) 553--596.

\bibitem {BGN2} C. Bai, L. Guo, X. Ni, $\O$-operators on associative algebras and associative Yang-Baxter equations,
Pacific J. Math. 256 (2012) 257--289.


\bibitem {B2} C. Bai, A unified algebraic approach to the classical Yang-Baxter equation, J. Phys. A 40 (2007)
11073--11082.


\bibitem{BD} A.A. Belavin, V.G. Drinfel'd, Solutions of the classical Yang-Baxter equation for simple Lie algebras, Funct. Anal. Appl. 16 (3) (1983) 159--180.

%\bibitem {BGP} P. Benito, V. Gubarev, A. Pozhidaev, Rota-Baxter operators on quadratic algebras, Mediterr. J. Math. 15 (2018) 18923 pp.

\bibitem{BPZ} A.A. Belavin, A.M. Polyakov, A.B. Zamolodchikov, Infinite conformal symmetry in two-dimen-sional
quantum field theory, Nucl. Phys. 241 (1984) 333--380.


\bibitem{BD2} A.A. Beilinson, V.G. Drinfeld, Chiral algebras. In: Amer. Math. Soc. Colloquium Publications,
vol. 51. AMS, Providence, RI, 2004.

\bibitem {BDK} B. Bakalov,  A. D'Andrea, V.G. Kac,  Theory of finite pseudoalgebras, Adv. Math. 162(1) (2001) 1--140.



\bibitem {BKV} B. Bakalov, V.G. Kac, A. Voronov, Cohomology of conformal algebras, Commun. Math. Phys. 200 (1999) 561--598.


\bibitem {BKL1}  C. Boyallian, V.G. Kac, J.I. Liberati, On the classification of subalgebras of $Cend_N$ and $gc_N$, J.
Algebra 260 (1) (2003) 32--63.

\bibitem {BKL2}  C. Boyallian, V.G. Kac, J.I. Liberati, Finite growth representations of infinite Lie conformal algebras,
J. Math. Phys. 44 (2) (2003) 754--770.


\bibitem {BO} R.E. Borcherds, Vertex algebras, Kac-Moody algebras, and the Monster, Proc. Nat. Acad. Sci. USA
83 (1986) 3068--3071.

\bibitem {CK} S.-J. Cheng, V.G. Kac, Conformal modules, Asian
J. Math. 1(1) (1997) 181--193. Erratum: Asian J. Math. 2 (1) (1998) 153--156.

%\bibitem{CKW1} S.-J. Cheng,  V. Kac, M. Wakimoto, Extensions of conformal modules, {\it In Topological Field Theory, Primitive Forms and Related Topics, Progr. Math.} Vol. 160., Kashiwara, M., et al. Eds. Boston: Birkh\"{a}user, 1998; 79-129.

%\bibitem {CKW2} S.-J. Cheng,  V. Kac, M. Wakimoto, Extensions of Neveu-Schwarz conformal modules, {\it
%J. Math. Phys.}, {\bf 41} (4) (2000) 2271--2294.


\bibitem{CGM} J.F. Carinena, J. Grabowski, G. Marmo, Quantum Bi-Hamiltonian systems, Int. J.
Mod. Phys. A 15 (2000) 4797--4810.

\bibitem{C} P. Cartier, On the structure of free Baxter algebras, Adv. Math. 9 (1972) 253--265.

\bibitem{DK} A. D'Andrea, V. Kac, Structure theory of finite
conformal algebras, Sel. Math., New Ser. 4 (1998) 377--418.

\bibitem{Das} A. Das, Deformations of associative Rota-Baxter operators, J. Algebra 560 (2020) 144--180.


\bibitem{Dor} I. Ya. Dorfman,  Deformations of Hamiltonian structures and integrable systems, In: Nonlinear and Turbulent Processes in Physics, ed. V. E. Zakharov, pp. 1313--1318. Gordon Breach, New York, 1984.


\bibitem{D} I.A. Dolguntseva,  The Hochschild cohomology for associative conformal algebras, Algebra
Logic 46 (2007) 373--384.

\bibitem{D2} I.A. Dolguntseva, Triviality of the second cohomology group of the conformal algebras
${\rm Cend}_n$ and ${\rm Cur}_n$, St. Petersburg Math. J. 21(1) (2009) 53--63.

\bibitem{DK2}  A. De Sole, V.G. Kac, Subalgebras of $gc_N$ and Jacobi polynomials, Canad. Math. Bull. 45 (4)
(2002) 567--605.


%\bibitem{EMP} K. Ebrahimi-Fard, D. Manchon, F. Patras, New identities in dendriform algebras, J. Algebra 320 (2008) 708--727.

%\bibitem{EM} K. Ebrahimi-Fard D. Manchon, Dendriform equations, J. Algebra 322 (2009) 4053--4079.


%\bibitem{F1} L. Foissy, Bidendriform bialgebras, trees, and free quasi-symmetric function, J. Pure Appl. Algebra 209 (2)
%(2007) 439--459.

\bibitem{F2}  I.B. Frenkel, J. Lepowsky, A. Meurman, Vertex operator algebras and the Monster, Pure and Applied
Mathematics, vol. 134, Academic Press, New York, 1998.


\bibitem{GD1} I.M. Gel'fand, I.Ya. Dorfman, Hamiltonian operators and algebraic structures related to them, Funct. Anal. Appl. 13 (4) (1979) 248--262.

\bibitem{GD2} I.M. Gel'fand, I.Ya. Dorfman, Schouten bracket and Hamiltonian operators, Funct. Anal. Appl. 14 (3)(1980) 71--74.


%\bibitem{GK1} V. Gubarev, P. Kolesnikov,  Embedding of dendriform algebras into Rota-Baxter algebras. Cent.
%Eur. J. Math. 11 (2) (2013) 226--245.

%\bibitem{GK2} Guo, L., Keigher, W. (2000). Baxter algebras and shuffle products. Adv. Math. 150(1):117--149.
\bibitem{G} L. Guo,  An Introduction to Rota-Baxter Algebra, Surveys of Modern Mathematics, Vol. 4.
Somerville: Intern. Press, 226 pp, 2012.

\bibitem{Ger} M. Gerstenhaber, The cohomology structure of an associative ring, Ann. of Math. 78 (1963) 267--288.


\bibitem{HB1} Y.Y. Hong,  C.M. Bai, On antisymmetric infinitesimal conformal bialgebras, J Algebra 586 (2021) 325--356.

\bibitem{HB2} Y.Y. Hong,  C.M. Bai, Conformal classical Yang-Baxter equation, $S$-equation and $\mathcal{O}$-operators, Lett. Math. Phys. 110 (2020) 885--909.

\bibitem{H1} Y.Y. Hong, Extending structures for associative conformal algebras, Linear Multilinear Algebra 67 (2019) 196--212.



\bibitem{KK} P.S. Kolesnikov, R.A. Kozlov, On the Hochschild cohomologies of associative
conformal algebras with a finite faithful representation, Commun. Math. Phys. 369 (2019) 351--370.

\bibitem{KAC} V.G. Kac, Vertex algebras for beginners, 2nd edn. In: University Lecture Series, vol. 10. AMS,
Providence, 1996.

\bibitem{KAC2} V.G. Kac, Formal distribution algebras and conformal algebras, XIIth International Congress in
Mathematical Physics (ICMP'97) (Brisbane), International Press, Cambridge, MA, 1999, pp. 80--97.

\bibitem{KO} P.S. Kolesnikov, Associative conformal algebras with finite faithful
representation, Adv. Math. 202 (2006) 602--637.

\bibitem{KS1} Y. Kosmann-Schwarzbach, From Poisson algebras to Gerstenhaber algebras, Ann. Inst. Fourier (Grenoble) 46(5) (1996) 1243--1274.

\bibitem{KS2} Y. Kosmann-Schwarzbach, Derived brackets, Lett. Math. Phys. 69 (2004) 61--87.

\bibitem{Ku} B.A. Kupershmidt, What a classical $r$-matrix really is? J. Nonlinear Math. Phy. 6 (1999) 448--488.


\bibitem{Loday}  J.-L. Loday, Dialgebras, in: Dialgebras and related operads, Lecture Notes in
Math. 1763 (2002) 7--66.

\bibitem{L} P. Leroux,  Construction of Nijenhuis operators and dendriform trialgebras, Int. J.
Math. Math. Sci. 49--52 (2004) 2595--2615.

%\bibitem{Lam} C.S. Lam, Decomposition of time-ordered products and path-ordered exponentials, J. Math. Phys. 39 (1998)
%5543--5558.

%\bibitem{LHB} Li, X. X., Hou, D. P., Bai, C. M. (2007). Rota-Baxter operators on pre-Lie algebras. J. Nonlinear Math.
%Phys. 14(2):269--289.

\bibitem{Li} J. Liberati, On conformal bialgebras, J. Algebra 319 (2008) 2295---2318.


\bibitem{LBS} J.F. Liu, C.M. Bai, Y.H. Sheng, Compatible $\O$-operators on bimodules over
associative algebras, J. Algebra 532 (2019) 80--118.


\bibitem{Ro} G.-C. Rota,  Baxter algebras and combinatorial identities I, Bull. Am. Math. Soc. 75 (2) (1969) 325--329.


\bibitem{Re1} A. Retakh, Associative conformal algebras of linear growth, J. Algebra 237 (2) (2001) 769--788.
\bibitem{Re2} A. Retakh, On associative conformal algebras of linear growth II,  J. Algebra 304 (1) (2006) 543--556.




\bibitem{U1} K. Uchino, Quantum analogy of Poisson geometry, related
Dendriform algebras and Rota-Baxter operators, Lett. Math. Phys. 85 (2008) 91--109.

\bibitem{U2} K. Uchino, Twisting on associative algebras and Rota-Baxter type operators, J. Noncommut. Geom. 4 (2010) 349--379.


\bibitem{U3} K. Uchino, Derived bracket construction and Manin products, Lett. Math. Phys. 93(1) (2010) 37--53.

%\bibitem{V1} Th. Voronov, Higher derived brackets and homotopy algebras, J. Pure Appl. Algebra 202(1-3) (2005) 133--153.



\bibitem{STS} M.A. Semenov-Tyan-Shanskii,  What is a classical $R$-matrix? Funktsional. Anal. i Prilozhen 17(4) (1983) 17--33.



\bibitem{SW} P. $\check{S}$evera, A. Weinstein, Poisson geometry with a 3-form background, Prog. Theor. Phys. Suppl. 144 (2001) 145--154



\bibitem{Wu} Z.X. Wu, Lie algebra structures oncohomology complexes of some $H$-pseudoalgebras, J. Algebra 396 (2013) 117--142.



\bibitem{YZ} L.M. Yuan, J.F. Liu, $\mathcal{O}$-operators, twisted Rota-Baxter operators and Nijenhuis operators on Lie conformal algebras, preprint, 2022.


\bibitem{Z1} E.I. Zelmanov, On the structure of conformal algebras, International Conference on Combinatorial
and Computational Algebra, Hong Kong, May 24--29, 1999, Contemp. Math. 264 (2000) 139--153.

\bibitem{Z2} E.I. Zelmanov, Idempotents in conformal algebras, Proceedings of the Third International Algebra
Conference, in: Y. Fong, et al. (Eds.),  pp. 257--266, 2003.


\bibitem{ZGG} T. Zhang, X. Gao,  L. Guo, Reynolds operators and their free objects from bracketed words and rootes trees, J. Pure Appl. Algebra 225 (12) (2021) 106766.




%%%%%%%%%%%%%%%%%%%%%%%%%%%%%%%%%%%%%%%%%%%%%%%%%%%%%%%%%%%%%%%%%%%%%%%%%%%%%%%%%%%%%%%
%%%%%%%%%%%%%%%%%%%%%%%%%%%%%%%%%%%%%%%%%%%%%%%%%%%%%%%%%%%%%%%%%%%%%%%%%%%%%%%%%%%%%%%%%%%%

%\bibitem{BBGW} Bai CM, Bai RB, Guo L, Wu Y. Transposed Poisson algebras, Novikov-Possion algebras and $3$-Lie algebras. arXiv:2005.01110v1.


%\bibitem{BN} Balinskii AA, Novikov SP. Poisson brackets of
%hydrodynamic type, Frobenius algebras and Lie algebras. Sov
%Math Dokl. 1985;32:228--231.



%\bibitem{BLLM} J. Bell, S. Launois, O. Le¨®n S¨¢nchez, R. Moosa, Poisson algebras via model theory and differential algebraic geometry, {\it J. Eur. Math. Soc.} {\bf19} (2017) 2019--2049.


%\bibitem{CK} Cantarini N, Kac V. Classification of linearly compact simple Jordan and generalized Poisson
%superalgebras. J Algebra. 2007;313:100--124.
%\bibitem{Ger}Gerstenhaber M.  The cohomology structure of an associative ring. Ann. Math. 1963;78:267--288.


%\bibitem{GK} Ginzburg V, Kaledin D. Poisson deformations of symplectic quotient singularities. Adv Math. 2004;186:1--57.

%\bibitem{H} Huebschmann J. Poisson cohomology and quantization. J Reine Angew Math. 1990;408:57--113.
%\bibitem{HLS} Hartwig JT, Larsson D, Silvestrov SD.
    %Deformations of Lie algebras using $\sigma$-derivation. J
            %Algebra. 2006;295:314--361.

%\bibitem{HWX} Han X, Wang DY, Xia CG. Linear commuting maps and biderivations on the Lie algebras $W(a,b)$. J Lie Theory. 2016;26:777--786.

%\bibitem{J}Jacobson N. Lie algebras. New York (NY): Interscience Publishers; 1962.

%\bibitem{Kay1} Kaygorodov IB. $\delta$-derivations of simple finite-dimensional Jordan superalgebras. Algebra Logic. 2007;46:318--329.
%\bibitem{Kay2} Kaygorodov IB. $\delta$-derivations of classical Lie superalgebras. Sib Math J. 2009;50:434--449.


%\bibitem{Kay3}Kaygorodov I, Shestakov I, Umirbaev U. Free generic Poisson fields and algebras. Comm Algebra. 2018;46:1799--1812.


%\bibitem{Kay4}Kaygorodova  I, Khrypchenko M. Poisson structures on finitary incidence algebras. J Algebra. 2021;578:402--420.


%\bibitem{K} Kontsevich M. Deformation quantization of Poisson manifolds. Lett Math Phys. 2003;66:157--216.


%\bibitem{KS}Kosmann-Schwarzbach Y. From Poisson to Gerstenhaber algebras. Ann Inst Fourier. 1996;46:1243--1274.


%\bibitem{Li} Li LC. Classical $r$-matrices and compatible Poisson structures for Lax equations on Poisson alge-bras. Commun Math Phys. 1999;203:573--592.


%\bibitem{MR} Markl M, Remm E. Algebras with one operation including Poisson and other Lie-admissible algebras. J Algebra. 2006;299:171--189.



%



%\bibitem{Xu1}Xu XP. On simple Novikov algebras and their
%irreducible modules. J Algebra. 1996;185:905--934.


%\bibitem{Xu2}Xu XP. Novikov-Poisson algebras. J Algebra. 1997;190:253--279.



%\bibitem{YYZ} Yao Y, Ye Y, Zhang P. Quiver Poisson algebras. J Algebra. 2007;312:570--589.

%\bibitem{Ze} Zelmanov E. On a class of local translation invariant Lie
%algebras. Sov Math Dokl. 1987;35:216--218.








}
\end{thebibliography}
\end{document}